\documentclass{gtart_h}

\def\ifplaintex{\expandafter\ifx\csname documentclass\endcsname\relax}

\def\gtp{{\mathsurround=0pt\it $\cal G\mskip-2mu$eometry \&\ 
$\cal T\!\!$opology $\cal P\!$ublications}}  

\def\recd{{\small Received:\qua\receiveddate\ifx\reviseddate\relax
\else\qquad Revised:\qua\reviseddate\fi\par}} 

\def\lognumber#1{\def\thelognumber{#1}}
\def\volumenumber#1{\def\thevolumenumber{#1}}
\def\volumeyear#1{\def\thevolumeyear{#1}}
\def\papernumber#1{\def\thepapernumber{#1}}
\def\pagenumbers#1#2{\def\startpage{#1}\def\finishpage{#2}}
\def\published#1{\def\publishdate{#1}}

\def\received#1{\def\receiveddate{#1}}

\def\accepted#1{\def\accepteddate{#1}}

\def\asciiaddress#1{\def\theasciiaddress{#1}}
\def\asciiemail#1{\def\theasciiemail{#1}}

\long\def\asciiabstract#1{\long\def\theasciiabstract{#1}}

\let\\\par\let\thelognumber\relax\let\thevolumenumber\relax
\let\thepapernumber\relax\let\thevolumeyear\relax\let\startpage\relax
\let\finishpage\relax\let\publishdate\relax\let\receiveddate\relax
\let\reviseddate\relax\let\accepteddate\relax\let\theasciititle\relax
\let\theasciiauthors\relax\let\theasciiaddress\relax
\let\theasciiabstract\relax

\let\theasciiemail\relax

\ifplaintex
\font\logobig=cmssbx10 scaled 3836
\font\logomed=cmssbx10 scaled 2557
\else
\font\logobig=cmssbx10 scaled 4200
\font\logomed=cmssbx10 scaled 2800
\fi

\long\def\makeagttitle{   
\count0=\startpage
\agt\hfill      
\hbox to 45truept{\vbox to 0pt{\vglue -13truept{\logomed A\kern -.37em{\logobig 
T}\kern -.38em G}\vss}\hss}
\break
{\small Volume \thevolumenumber\ (\thevolumeyear)
\startpage--\finishpage\nl
Published: \publishdate}

\vglue .25truein

{\parskip=0pt\leftskip 0pt plus
1fil\def\\{\par\smallskip}{\Large\bf\thetitle}\par\medskip} \vglue
0.05truein

{\parskip=0pt\leftskip 0pt plus 1fil\def\\{\par}{\sc\theauthors}
\par\medskip}%
 
\vglue 0.03truein 

{\small\leftskip 25truept\rightskip 25truept{\bf Abstract}\stdspace\theabstract

{\bf AMS Classification}\stdspace\theprimaryclass
\ifx\thesecondaryclass\relax\else; \thesecondaryclass\fi\par
{\bf Keywords}\stdspace \thekeywords\par}\vglue 7truept

}   

\ifplaintex
\font\phead=cmsl9 scaled 950
\font\pnum=cmbx10 scaled 913
\font\pfoot=cmsl9 scaled 950
\headline{\vbox to 0pt{\vskip -4.5mm\line{\small\phead\ifnum
\count0=\startpage ISSN 1472-2739 (on-line) 1472-2747 (printed)
\hfill {\pnum\folio}\else\ifodd\count0\def\\{ }%
\ifx\theshorttitle\relax\thetitle\else\theshorttitle\fi\hfill{\pnum\folio}
\else\def\\{ and }{\pnum\folio}\hfill\ifx\theshortauthors\relax\theauthors
\else\theshortauthors\fi\fi\fi}\vss}}
\footline{\vbox to 0pt{\vglue 0mm\line{\small\pfoot\ifnum\count0=\startpage
\copyright\ \gtp\hfill\else
\agt, Volume \thevolumenumber\ (\thevolumeyear)\hfill\fi}\vss}}
\else
\font=cmsl9 scaled 1050
\font\lnum=cmbx10 
\font=cmsl9 scaled 1050
\makeatletter
\def\@oddhead{{\small\lhead\ifnum\count0=\startpage ISSN 1472-2739 
(on-line) 1472-2747 (printed)\hfill {\lnum\number\count0}\else\ifodd\count0
\def\\{ }\ifx\theshorttitle\relax \thetitle \else\theshorttitle\fi\hfill
{\lnum\number\count0}\else\def\\{ and }{\lnum\number\count0}
\hfill\ifx\theshortauthors\relax 
\theauthors\else\theshortauthors\fi\fi\fi}}\def\@evenhead{\@oddhead}
\def\@oddfoot{\small\lfoot\ifnum\count0=\startpage\copyright\ \gtp\hfill\else
\agt, Volume \thevolumenumber\ (\thevolumeyear)\hfill\fi}
\def\@evenfoot{\@oddfoot}
\makeatother
\fi
\let\maketitlepage\makeagttitle
\let\makeshorttitle\maketitlepage
\let\maketitle\maketitlepage

\newwrite\gtoutfile
\long\gdef\makeheadfile{  
{\def\\{, }\def\s{ }
\immediate\openout\gtoutfile head.xxx
\immediate\write\gtoutfile{Proxy-for: \ifx\theasciiauthors\relax
\theauthors\else\theasciiauthors\fi\s<\ifx\theasciiemail\relax\theemail\else\theasciiemail\fi>}
\immediate\write\gtoutfile{\noexpand\\}
\immediate\write\gtoutfile{Authors: \ifx\theasciiauthors\relax
\theauthors\else\theasciiauthors\fi}
{\def\\{ }\immediate\write\gtoutfile{Title: \ifx\theasciititle\relax
\thetitle\else\theasciititle\fi}}
\immediate\write\gtoutfile{Subj-class: GT or SG, GR etc}
\immediate\write\gtoutfile{MSC-class: \theprimaryclass\ifx\thesecondaryclass\relax\else, \thesecondaryclass\fi}
\immediate\write\gtoutfile{Journal-ref: Algebr. Geom. Topol. \thevolumenumber\s
(\thevolumeyear) \startpage-\finishpage}
\immediate\write\gtoutfile{Comments: Published by Algebraic and
Geometric Topology at}
\immediate\write\gtoutfile{\s\s\s  http://www.maths.warwick.ac.uk/agt/AGTVol\thevolumenumber/agt-\thevolumenumber-\thepapernumber.abs.html}
\immediate\write\gtoutfile{\noexpand\\}
\immediate\write\gtoutfile{}
\ifx\theasciiabstract\relax
\immediate\write\gtoutfile{\theabstract}\else
\immediate\write\gtoutfile{\theasciiabstract}\fi
\immediate\write\gtoutfile{}
\immediate\write\gtoutfile{\noexpand\\}
\immediate\write\gtoutfile{}
\immediate\closeout\gtoutfile}}  

\def\maketitlepage{\makeagttitle\makeheadfile}
\let\makeshorttitle\maketitlepage
\let\maketitle\maketitlepage

\lognumber{35}
\volumenumber{5}
\volumeyear{2005}
\papernumber{35}
\pagenumbers{835}{864}
\received{20 October 2004} 
\accepted{13 June 2005}
\published{24 July 2005}

\usepackage{amssymb, amsmath, epsf, psfrag}

\def\psfraga <#1,#2> #3#4{%
\psfrag {#3}{\smash{\rlap{\kern #1 \raise #2\hbox{#4}}}}}

\newtheorem{pro}{Proposition}[section]
\newtheorem{thm}[pro]{Theorem}
\newtheorem{lem}[pro]{Lemma}

\theoremstyle{definition}

\theoremstyle{remark}

\def\injrad{{\rm injrad}}
\def\Diameter{{\rm Diameter}}
\def\length{{\rm length}}
\def\Core{{\rm Core}}
\def\Hull{{\rm Hull}}
\def\Isom{{\rm Isom}}

\begin{document}

\title{Counting immersed surfaces\\in hyperbolic 3-manifolds}

\authors{Joseph D. Masters}
\address{Mathematics Department, Rice University, Houston TX 77005, 
USA\\{\rm Current address:}\\
Mathematics Department, SUNY Buffalo, Buffalo NY 14260, USA}

\asciiaddress{Mathematics Department, Rice University, Houston TX 77005, USA\\
Current address:\\
Mathematics Department, SUNY Buffalo, Buffalo NY 14260, USA}

\gtemail{\mailto{mastersj@rice.edu}, \mailto{jdmaster@buffalo.edu}}
\asciiemail{mastersj@rice.edu, jdmaster@buffalo.edu}

\begin{abstract}
 We count the number of conjugacy classes of maximal, genus $g$,
 surface subroups in hyperbolic 3-manifold groups.  For any closed
 hyperbolic 3-manifold, we show that there is an upper bound on this
 number which grows factorially with $g$.  We also give a class of
 closed hyperbolic 3-manifolds for which there is a lower bound of the
 same type.
\end{abstract}

\asciiabstract{%
We count the number of conjugacy classes of maximal, genus g,
surface subroups in hyperbolic 3-manifold groups.  For any closed
hyperbolic 3-manifold, we show that there is an upper bound on this
number which grows factorially with g.  We also give a class of
closed hyperbolic 3-manifolds for which there is a lower bound of the
same type.}

\primaryclass{57M50}\secondaryclass{57N16, 57M27}
\keywords{Surface subgroups, bending, pleated surfaces, reflection orbifolds}

\maketitle

\section{Introduction}
 A major problem in the study of hyperbolic 3-manifolds is to determine
 the extent to which they contain useful surfaces.  It is conjectured
 that every closed hyperbolic 3-manifold contains a number
 of immersed incompressible surfaces, each of which lifts
 to an embedding in a finite sheeted cover.
 The aim of this paper is to count the number in some appropriate sense.
 
 We fix attention on a closed hyperbolic 3-manifold (or orbifold) $M$.
 Since we are interested in immersed surfaces up to homotopy,
 we shall be counting surface subgroups of $\pi_1 M$ only up to conjugacy.
  Moreover, given our topological motivation,
 it is natural not to distinguish between a given immersed surface
 and one which covers it, and this leads to different
 ways of counting.

 We define $s(M,g)$ to be the number of conjugacy classes of 
 \textit{maximal}  surface groups
 of genus at most $g$ in $\pi_1(M)$.
  Two subgroups  of $\pi_1 M$ are said to be \textit{commensurable}
 if their intersection has finite index in both.
 Let $s_1(M,g)$ denote the number of classes of surface
 groups in $\pi_1 M$ under the equivalence relation
 generated by conjugacy and commensurability.
 Let $s_2(M,g)$ be the total number of conjugacy classes of surface
 subgroups in $\pi_1 M$ of genus at most $g$.

 Thurston proved (Corollary 8.8.6 of \cite{T}) that $s_2(M,g)$ is
 finite, but the proof does not give a practical upper bound.
 In \cite{S}, Soma gave an explicit upper bound, which is doubly
 exponential in the genus.  We show:

\begin{thm} \label{upper}
Let $M$ be a closed hyperbolic 3-manifold.  Then there exists a constant
 $c_2 > 0$ such that, for large $g$,
 \begin{equation*}
s_2(M,g) < e^{c_2 g \log g}.
\end{equation*}
\end{thm}

 We say that a totally geodesic immersion $f:S \rightarrow M$,
 is \textit{transverse} if its image has a curve of transverse
 self-intersection; thus $f$ is not a finite cover of an embedding.

\begin{thm} \label{main}
Let $M$ be a closed hyperbolic 3-manifold, and suppose there
 is a transverse, totally geodesic immersion $f:S \rightarrow M$,
 for some hyperbolic surface $S$.
 Then there exist constants $c_1, c_2 > 0$ such that, for large $g$,
\begin{equation*} 
e^{c_1 g \log g} < s_1(M,g) \leq s(M,g) \leq s_2(g,M) <
                     e^{c_2 g \log g}.
\end{equation*}
\end{thm}

 Theorem \ref{main} may be viewed as a co-dimension one analogue of 
 the Prime Geodesic Theorem for closed hyperbolic 3-manifolds
 (see \cite{EGM}), which states that the number
 of conjugacy classes of maximal 1-manifold groups grows exponentially
 with the length of the corresponding geodesics.
 The 3-manifolds of Theorem \ref{main}, in this sense,
 grow faster in two dimensions than in one dimension.
 It is also interesting to count \textit{finite-index} subgroups
 of $\pi_1 M$-- see \cite{Lu}.  For results on counting
 \textit{totally geodesic} immersions of surfaces,
 see \cite{MR} and \cite{JM}.  Anneke Bart and Brian Mangum
 have also obtained results on cutting and pasting immersed surfaces.

\medskip
\textbf{Organization}\qua
 In Section 2 we prove the upper bound. 
 As in Thurston's finiteness proof, the first step is
 to homotop each immersed surface into pleated form.
 The rest of the proof relies on finding
 triangulations for the surfaces whose edges are sufficiently
 short, and counting the possible graphs which arise as 1-skeleta.

 The proof of the lower bound is by an explicit construction
 of immersed incompressible surfaces, which is related to Thurston's bending
 deformation.
 We collect some preliminary material on surfaces and graphs in Section 3,
 and give the construction in Section 4.
 The main difficulty is to show that the resulting surface groups are
 inequivalent in the sense defined above.  We do this by proving
 that they correspond to hyperbolic manifolds with non-isometric convex cores.

 Section 5 is devoted to examples where
 $M$ is a right-angled reflection orbifold.
 In this case, we are able to give more explicit, combinatorial proofs
 of both bounds.  In particular, the proof of the upper bound
 avoids the usual theory of pleated surfaces, and as a corollary, we get
 an alternate proof of Thurston's finiteness theorem for these manifolds.

\begin{thm} \label{reforbs}
Let $P \subset \mathbb{H}^3$ be a compact,
 right-angled polyhedron, and let $M$ be the associated
 reflection orbifold.
 Then for the constants of Theorem \ref{main},
 we may take $c_1 = 1$ and $c_2 = 8c(P)+1$,
 where $c(P)$ is the maximum number of edges in a face disk of $P$
 (see Section 5 for definitions).
\end{thm}
 
Thanks are due to Peter Shalen for a useful conversation,
 to Lewis Bowen for pointing out an error in a previous version,
 and to the University of Texas at Austin for its hospitality
 during the completion of this work.

\section{An upper bound}

In this section, we prove Theorem \ref{upper}.

If $G$ is a graph, let $\mathcal{V}(G)$ denote
 the vertex set of $G$, and let $\mathcal{E}(G)$
 denote the edge set of $G$.

\begin{lem} \label{surface}
 There is a constant $k = k(s)$ such that
 any hyperbolic surface $S$ with $\injrad(S) \geq s$
 contains an embedded graph $G$ satisfying:

\begin{enumerate}
\item Every edge $e \in \mathcal{E}(G)$ is a geodesic arc of length $< s$,\\
\item $|\mathcal{V}(G)| < kg$,\\
\item $degree(v) < k$ for all vertices $v \in \mathcal{V}(G)$, and\\
\item $S-G$ is a disjoint union of open disks.
\end{enumerate}
\end{lem}

\begin{proof}
 Let $\mathcal{B}= \{ B_1, B_2, ..., B_n \}$ be a maximal collection of
 disjoint balls of radius $s/4$ in $S$, and let $v_i$ denote
 the center of ball $B_i$.
 If $x,y \in S$ are points with $d(x,y) < s$, then there is a unique
 shortest arc in $S$ connecting them, which we
 denote $e(x,y)$.  

Let
\begin{eqnarray*}
 \mathcal{S} &=& \{ (v_i, v_j) | int \, B(v_i,s/2) \cap int \, B(v_j,s/2)
                       \neq \emptyset \},\\
 \mathcal{E} &=& \{ e(v_i,v_j) | (v_i,v_j) \in \mathcal{S} \}
\end{eqnarray*}
 We define a subset $G$ of $S$ by:
\begin{eqnarray*}
G = \bigcup_{e \in \mathcal{E}}  e,
\end{eqnarray*}
and we give $G$ the structure of a graph by declaring the vertex set
 of $G$ to be:
\begin{eqnarray*}
\mathcal{V}(G) = \bigcup_{e,e^{\prime} \in \mathcal{E}}
  e \cap e^{\prime},
\end{eqnarray*}
and then defining the edges of $G$ to be closures of components
 of $G - \mathcal{V}(G)$.

\begin{lem}
There is a constant $k = k(s)$, (independent
 of $\mathcal{B}$) such that $|\mathcal{V}(G)| < kg$,
 and $degree(v) < k$ for all $v \in \mathcal{V}(G)$.
\end{lem}

\begin{proof}
We have:
\begin{equation*}
 |\mathcal{B}| \leq \frac{Area(S)}{Area(B(x,s/4))}
 = \frac{4\pi(g-1)}{4\pi sinh^2(\frac{s}{8})},
\end{equation*}
so $|\mathcal{B}| < m_1g$ for some $m_1$.
 If $B \subset S$ is a fixed $s/4$-ball,
 there are at most a fixed number-- call it  $m_2$--
 of disjoint $s/4$-balls which will fit within an $s$-neighborhood of $B$.
 So we have
\begin{eqnarray*}
|\mathcal{E}| = |\mathcal{S}| < m_2|\mathcal{B}| < m_1m_2 g.
\end{eqnarray*}

 Let $e \in \mathcal{E}$. Since $\length(e) < s$,
 there are at most a fixed number-- call it  $m_3$-- of
 disjoint $s/4$-balls in $S$ whose centers
 can fit within an $s$-neighborhood of $e$.
 Therefore, $e$ can intersect at most
 $\binom{m_3}{2}$ different $e(v_i, v_j)$'s.
 Since no pair of geodesic arcs of length less than $s$ can
 intersect more than once in $S$, we have
 \begin{eqnarray*}
|\mathcal{V}(G)| &=& |\bigcup_{e,e^{\prime} \in \mathcal{E}}
 e \cap e^{\prime}|,\\
 &\leq& \binom{m_3}{2} |\mathcal{E}|\\
 &<& m_1m_2 \binom{m_3}{2} g.
\end{eqnarray*}
The fact that $e$ can intersect at most
 $\binom{m_3}{2}$ other edges of $\mathcal{E}$ also implies
\begin{eqnarray*}
deg(v) < 2\binom{m_3}{2} + 2,
\end{eqnarray*}
for all $v \in \mathcal{V}(G)$.
 Therefore we may take $k = Max\{m_1m_2\binom{m_3}{2}, 2\binom{m_3}{2} +2\}$.
\end{proof}

\begin{lem}
The surface $S-G$ is a disjoint union of open disks.
\end{lem}

\begin{proof}
 Let $\mathcal{C} = \{B(v_1, s/2), ..., B(v_n, s/2) \}$. 
 Since $\mathcal{B}$ is maximal, $\mathcal{C}$ covers $S$, and therefore
 it is enough to show that every ball $B(v_i, s/2) \in \mathcal{C}$ is covered
 by (the closures of) simply connected components of $S-G$.
 
 The maximality of $\mathcal{C}$ implies that $\mathcal{C} - B(v_i,
 s/2)$ covers $\partial B(v_i,s/2)$.  Let $B_{i1}, ..., B_{im} \in
 \mathcal{C}-B(v_i,s/2)$ be balls which cover $\partial B(v_i,s/2)$,
 arranged such that $B_{ij} \cap B_{i(j+1)} \neq \emptyset$ for all
 $j$.  Then each triple $(v_i,v_{ij},v_{i(j+1)})$ is the vertex set of
 a geodesic triangle, $T_j$, contained in an $s$-ball (embedded, since
 $\injrad(S) \geq s$).  Thus each $T_j$ bounds a disk $D_j$, which is
 the union of (closures of) components of $S-G$. The ball $B(v_i,s/2)$
 is covered by $\{D_1, ..., D_m\}$.
\end{proof}

This concludes the proof of Lemma \ref{surface}.
\end{proof}

 We now return to the proof of the upper bound.
By \cite{T}, given any $\pi_1$-injective immersion $f: S \rightarrow M$,
 we can find a hyperbolic structure on $S$, and a homotopy of $f$ so that
 it is pleated with respect to this structure.  
 Let $s = \injrad(M)$; since $f$ is pleated
 it takes closed loops in $S$ to closed loops of
 equal or shorter length in $M$, and therefore $\injrad(S) \geq s$.
 There is then a graph $G \subset S$,
 with the properties stated in Lemma \ref{surface}.

 Suppose $f^{\prime}:S^{\prime} \rightarrow M$ is another
 $\pi_1$-injective, pleated map of a genus $g$ surface,
 with a graph $G^{\prime} \subset S^{\prime}$
 given by Lemma \ref{surface}.  With slight abuse of terminology,
 we say that $f$ and $f^{\prime}$ are \textit{homotopic} if
 there is a map $g: S^{\prime} \rightarrow S$
 such that $fg$ is homotopic to $f^{\prime}$.

 Suppose there is a bijection
 $h: \mathcal{V}(G) \rightarrow \mathcal{V}(G^{\prime})$,
 and label the vertices $\mathcal{V}(G) = \{v_1, ..., v_n \}$,
 and $\mathcal{V}(G^{\prime}) = \{v_1^{\prime}, ..., v_n^{\prime} \}$,
 where $v_i^{\prime} = h v_i$.
 The map $h$ induces a graph $h(G)$ on $\mathcal{V}(G^{\prime})$
 by the rule that $v_i^{\prime}$ and $v_j^{\prime}$ are adjacent in
 $h(G)$ if and only if $v_i$ and $v_j$ are adjacent in $G$.

\begin{lem} \label{graphs}
Suppose that $d(f^{\prime} v_i^{\prime}, f v_i) < s/4$ for all i,
 and that $h(G) = G^{\prime}$.
 Then $f$ and $f^{\prime}$ are homotopic.
\end{lem}

\begin{proof}
Let $\delta_i$ be a segment of length less than $s/4$ connecting $f v_i$
 and $f^{\prime} v_i^{\prime}$.
  Sliding the $v_i^{\prime}$'s along the $\delta_i$'s,
 we may homotope $f^{\prime}|_{G^{\prime}}$
 to a map $g:G^{\prime} \rightarrow M$, so that
 $g(v_i^{\prime}) = f(v_i)$
 and so that $g G^{\prime}$ is contained in
 $f^{\prime} G^{\prime} \bigcup_i \delta_i$. 
 If $E^{\prime}$ is an edge of $G^{\prime}$, then its image under
 $g$ is given by:
\begin{equation*}
 g (E^{\prime}) = f^{\prime} (E^{\prime}) \cup \delta_i \cup \delta_j,
 \textnormal{  (for some i,j)},
\end{equation*}
 which has total length less than $s$.  Since $\injrad(M) = s$, any two
 segments of $M$ with the same endpoints which have length less than $s$ 
 are homotopic, fixing endpoints.
  Therefore, the map $f^{\prime}|_{G^{\prime}}$
 is homotopic to $fh|_{G^{\prime}}$.
 Since the complementary regions of $G$ and $G^{\prime}$
 are disks, the map $h$ can be extended to a map $h:S^{\prime} \rightarrow S$,
 such that $f^{\prime}$ and $fh$ are homotopic.
\end{proof}

 Let $\mathcal{V}=\{ v_1, ..., v_k \}$ be a fixed set of vertices, 
 and let $\mathcal{G}(\mathcal{V})$ be the set of all graphs on $\mathcal{V}$.
 Note that the relation of equality between graphs on $\mathcal{V}$
 is stronger than the relation of isomorphism between graphs on $\mathcal{V}$.

 Let 
\begin{eqnarray*}
\mathcal{G}(\mathcal{V},n) =
 \{ G \in \mathcal{G}(\mathcal{V}) | \textrm{ each vertex
 of } G \textrm{ has degree at most } n \}
\end{eqnarray*}
 A computation shows that
\begin{equation}
|\mathcal{G}(\mathcal{V}, n)| \leq |\mathcal{V}|^{n |\mathcal{V}|}.
\end{equation}
 
 Let $k_1=k(s)$ be the constant provided by Lemma \ref{surface}.
 Let $\mathcal{C}$ be a collection of balls of radius
 $s/4$ which covers $M$, and let $k_2 \geq |\mathcal{C}|$.
 Suppose we have a collection $\mathcal{S}$ of more than 
 $(k_1g)(k_1g)^{k_2}(k_1g)^{k_1^2g}$
 pleated maps $f_i:S_i \rightarrow M$, where $S_i$ is a hyperbolic
 surface of genus $g$.
 Associated to each $f_i \in \mathcal{S}$
 is a graph $G_i \subset S_i$ satisfying the properties
 stated in Lemma \ref{surface}. Since $|\mathcal{V}(G_i)| < k_1g$,
 there is a subset $\mathcal{S}_1 \subset \mathcal{S}$
 of $(k_1g)^{k_2} (k_1g)^{k_1^2 g}$ pleated maps whose graphs all
 have the same number of vertices.
 Since there are at most $(k_1g)^{k_2}$ ways to map
 $\mathcal{V}(G_i)$ to $\mathcal{C}$,
 then there is a subset $\mathcal{S}_2 \subset \mathcal{S}_1$ of
 $(k_1g)^{k_1^2 g}$ pleated maps, such that,
 if we fix a map $f \in \mathcal{S}_2$
 with associated graph $G$, then for any $f_i \in \mathcal{S}_2$,
 there are bijections $h_i: \mathcal{V}(G_i) \rightarrow \mathcal{V}(G)$,
 such that, for all $v \in \mathcal{V}(G_i), d(f h_i(v), f_i(v)) \leq s/4$.
 By Equation 1, there is a pair
 of distinct maps $f_i, f_j \in \mathcal{S}_2$
 such that $h_i(G_i) = h_j(G_j)$.
 Therefore, by Lemma \ref{graphs}, $f_i$ is homotopic to $f_j$.
 This concludes the proof of the upper bound.

\section{Properties of graphs and surfaces}

By a \textit{metric graph}, $G$, we mean a graph with a metric space
 structure, determined by the following procedure:
 we assign a length to each edge, $e$, and give
 $e$ a path metric induced from the interval $[0,\length(e)]$.
 The path metric on the edges then determines a path
 metric on $G$. The \textit{standard} metric
 on $G$ is obtained by setting the length
 of each edge to be one.
 A map $f:X \rightarrow Y$
 between metric spaces is a \textit{(k,c)-quasi-isometry} ($k,c>1$)
 if $\frac{d(fx_1,fx_2)}{k} - c < d(x_1,x_2) < kd(fx_1,fx_2)+c$
 for all $x_1,x_2 \in X$, and $d(y, f(X)) < c$ for all $y \in Y$.

\begin{lem} \label{quasi}
 Let $G$ and $G^{\prime}$ be metric graphs,
 with no vertices of degree 2, and with all
 edges longer than $u=6k^3(k^2+3)c$.
 Then if $G$ and $G^{\prime}$
 are  $(k,c)$-quasi-isometric, they are isomorphic.
\end{lem}

\begin{proof}
 Let $f: G \rightarrow G^{\prime}$
 be a $(k,c)$-quasi-isometry.  Then, there exist numbers
 $k^{\prime}, c^{\prime}$ and $s>0$, and a map
 $g: G^{\prime} \rightarrow G$,
 which is a $(k^{\prime},c^{\prime})$-quasi-isometry, such that
 $d(gfx,x) < s$, and $d(fgy,y) < s$
 for all $x \in G$ and $y \in G^{\prime}$.  In fact, we may take
 $k^{\prime} = k$, $c^{\prime} = 3kc$, and $s = kc$.

 Let $x$ be a vertex  of $G$ of degree $n > 2$,
 and let $e_1, ..., e_n$ be distinct edges incident to $x$.
 Since $f$ is a $(k,c)$-quasi-isometry, any point of $fe_i$ is within
 distance $c$ of another point of $fe_i$,
 and since $\length(e_i) > k(k^2+2)c$, then
 $\Diameter(fe_i) > k(k^2+2)c/k-c = (k^2+1)c$.
 Thus, for each $i$, there are points
 $x_i \in e_i$  such that
\begin{eqnarray*}
(k^2+1)c < d(fx_i, fx) < (k^2+1)c +c.
\end{eqnarray*}
 Then
\begin{eqnarray*}
d(x_i,x) &>& d(fx_i,fx)/k - \frac{c}{k}\\
&>& (k^2+1)c/k-\frac{c}{k}\\
&=& kc.
\end{eqnarray*}
 Since edges of $G$ all have length greater
 than $2d(x_i,x)$, then
\begin{eqnarray*}
 d(x_i,x_j)= d(x_i,x)+d(x_j,x) > 2kc,
\end{eqnarray*}
 $$d(fx_i,fx_j) > (2kc)/k - c > c.\leqno{\hbox{and so}}$$
 So, in the $(k^2+2)c$-neighborhood of $fx$,
 there are $n$ points $fx_i$, such that
 $|d(fx_i, fx) - d(fx_j,x)| < c$,
 and $d(fx_i, fx_j)> c$ for all $i \neq j$. 
 This implies that $fx$ is within distance $(k^2+2)c$ of a vertex
 of $G^{\prime}$.

 If $x$ is a degree 1 vertex, then
 we claim that $y=fx$ is
 within distance $k^{\prime}(k^{\prime 2}+2)c^{\prime}$
 of a degree 1 vertex of $G^{\prime}$.
 For suppose not. Then there is an interval $e$ of length
 $2k^{\prime}(k^{\prime 2}+2)c^{\prime}$
 isometrically embedded in $G^{\prime}$,
 such that $fx$ is the midpoint of $e$.
 Let $e_1$ and $e_2$ be sub-intervals of $e$ with
 $e_1 \cap e_2 = \{ fx \}$, and
 $\length(e_i) = k^{\prime}(k^{\prime 2}+2)c^{\prime}$.
 Then
 $\Diameter(ge_i) > [k^{\prime}(k^{\prime 2}+2)c^{\prime}]/k^{\prime}-c^{\prime}
       = (k^{\prime 2}+1)c^{\prime}$, so for both
 $i$'s there exists a point $y_i \in e_i$ such that
\begin{eqnarray*}
 (k^{\prime 2}+1)c^{\prime} < d(gy_i, gy)
   <   (k^{\prime 2}+1)c^{\prime} + c^{\prime}.
\end{eqnarray*}
 We have:
\begin{eqnarray*}
 d(y_i, y)&>& \frac{d(gy_i,gy)}{k^{\prime}} - \frac{c^{\prime}}{k^{\prime}}\\
 &>& \frac{(k^{\prime 2}+1)c^{\prime}}{k^{\prime}} - \frac{c^{\prime}}{k^{\prime}}\\
 &=& k^{\prime}c^{\prime}
\end{eqnarray*}
 Since all edges of $G^{\prime}$ have length greater
 than $d(y_i,y)$, then
\begin{eqnarray*}
  d(y_1, y_2) =d(y_1,y)+d(y_2,y) > 2k^{\prime}c^{\prime},
\end{eqnarray*}
$$d(gy_1,gy_2) > (2k^{\prime}c^{\prime})/k^{\prime}-c^{\prime}=c^{\prime}.
\leqno{\hbox{and}}$$
 However, we have $d(gy,x)=d(gfx,x)<s=kc$,
 so $gy$ and $x$ are contained in the same edge. Since $x$ has degree one,
 and
\begin{eqnarray*}
 d(gy_i,x)&>& d(gy_i,gy) - d(gy,x)\\
&>&(k^{\prime 2}+1)c^{\prime}-s\\
&=&(k^2+1)(3kc)-kc\\
&>&kc=s\\
&>&d(gy,x),
\end{eqnarray*}
 then $d(gy_i, x)=d(gy_i,gy)+d(gy,x)$.
 Then
\begin{eqnarray*}
 |d(gy_1,x)- d(gy_2,x)| &=& |d(gy_1,gy)+d(gy,x)-(d(gy_2,gy)+d(gy,x))|\\
&=&|d(gy_1,gy)-d(gy_2,y)|\\
&<& c^{\prime}.
\end{eqnarray*}
 Since
 $d(gy_i,x) =d(gy_i,gy)+d(gy,x)<(k^{\prime 2}+2)c^{\prime}+s=(k^2+2)(3kc)+kc$,
 then $gy_i$ and $x$ are contained in the same edge.
 Therefore
 $d(gy_1, gy_2) = |d(gy_1,x)-d(gy_2,x)|< c^{\prime}$, for a contradiction.
 Therefore, $fx$ is within distance $k^{\prime}(k^{\prime 2}+2)c^{\prime}$
  of a degree 1 vertex of $G^{\prime}$.

 Let
\begin{eqnarray*}
 t &=& Max( (k^2+2)c,k^{\prime}(k^{\prime 2}+2)c^{\prime})\\
  &=& Max( (k^2+2)c,k(k^2+2)3kc)\\
 &=&3k^2(k^2+2)c\\
 &<&u/2
\end{eqnarray*}
  We have shown $f(\mathcal{V}(G))
      \subset N_t(\mathcal{V}(G^{\prime}))$.
 Since $t<u/2$, there is an induced map
 $f^*:\mathcal{V}(G) \rightarrow \mathcal{V}(G^{\prime})$,
 defined by $f^*v = N_{u/2}(fv) \cap  \mathcal{V}(G^{\prime})$.
 Similarly, if
 $t^{\prime} = Max( (k^{\prime 2}+2)c,k(k^2+2)c)$, then
 $t^{\prime} <u/2$,
 $g(\mathcal{V}(G^{\prime}))
      \subset N_{t^{\prime}}(\mathcal{V}(G))$,
  and there is an induced map
 $g^*:\mathcal{V}(G^{\prime}) \rightarrow \mathcal{V}(G)$.

 We have:
\begin{eqnarray*}
d(gf^* v,v) &\leq& d(v, gfv)+ d(gfv, gf^*v)\\
 &\leq& s+k^{\prime}d(fv,f^*v)+c^{\prime}\\
 &\leq& s+k(3k^2(k^2+2)c)+3kc\\
 &<& u/2
\end{eqnarray*} 
 Therefore $g^*f^* v = v$, and similarly $f^* g^* v = v$
 for all $v \in \mathcal{V}(G^{\prime})$.
 So the maps $f^*$ and $g^*$
 are inverses to each other, and hence bijections.
  
 Furthermore, if $v_1,v_2$ are adjacent vertices in $G$,
 with an edge $e$ connecting them, then $f(e)$
 is a $(k,c)$-quasi-geodesic segment, intersecting
 $N_t(f^*v_1 \cup f^*v_2)$,  but disjoint
 from $N_t(\mathcal{V}(G^{\prime})-\{f^*v_1, f^*v_2 \})$,
 and this implies that $f^*v_1$ and $f^*v_2$ are adjacent
 in $G^{\prime}$.  A similar statement
 holds for $g^*$, and so the map $f^*$ induces
 an isomorphism from $G$ to $G^{\prime}$.
 \end{proof}

 Let $G$ be a graph,
 with a non-separating basepoint $p$ contained in the interior of an edge. 
 Let $N(p)$ be a regular open neighborhood of $p$ in $G$,
 and let $G/p = G - N(p)$.
 Then $G/p$ has the structure of a graph, with
 $\mathcal{V}(G/p) = \mathcal{V}(G) \cup \{v_1, v_2 \}$,
 where $v_1$ and $v_2$ are degree one vertices.
 Say that the pairs $(G, p)$ and
 $(G^{\prime}, p^{\prime})$ are \textit{inequivalent} if
 the universal cover of $G/p$ is  not isomorphic
 to the universal cover of $G^{\prime}/p^{\prime}$. 

\begin{lem} \label{graphinequiv}
Let $G$ be a bouquet of two circles,
 with basepoint $p$ in the interior of an edge.
  Then there is a constant $c>0$
 such that, for all $n$, the number
 of inequivalent covers $(\tilde{G},\tilde{p})$
 of $(G,p)$ of degree $n$ is at least $(cn)!$.
\end{lem}

\begin{proof} 
 Let $a,b$ be the edges of $G$, and assume $p \in b$.
 Any graph $\tilde{G}$ whose vertices all have degree four
 is a cover of $G$, and we may specify the covering map
 by directing the edges of $\tilde{G}$,
 and labeling each of them either $a$ or $b$.  Let $\sigma$ be
 an order two permutation on $\{ 1, ..., n \}$. We consider $4n$-fold covers
 $\pi: (\tilde{G},\tilde{p}) \rightarrow (G,p)$
 constructed as follows (see Figure \ref{fig:1}):
\begin{enumerate}
\item  Begin with $4n$ vertices $v_1, ..., v_{4n}$, and, for all $i$, connect
 $v_i$ and $v_{i+1}$ with an edge labeled $a$.
\item If either
\begin{itemize}
\item[(a)] $i$ is even
\item[(b)] $i \geq 2n$, or
\item[(c)]  $i = 2j-1$, with $j \leq n$, and $\sigma j = j$,
\end{itemize}
then attach a 1-cycle labeled $b$ to $v_i$.
\item If $\sigma(i) = j$, and $i \neq j$, then
 attach a 2-cycle with edges labeled $b$ to $v_{2i-1},v_{2j-1}$.
\item Choose a basepoint, $\tilde{p}$, in the 1-cycle
 which is adjacent to $v_0$ and labeled $b$.
\end{enumerate}

\begin{figure}[ht!]\small
\psfraga <-2pt, 0pt> {a}{$a$}
\psfraga <-2pt, -2pt> {b}{$b$}
\psfraga <0pt, 0pt> {p}{$p$}
{\epsfxsize=2.5in \centerline{\epsfbox{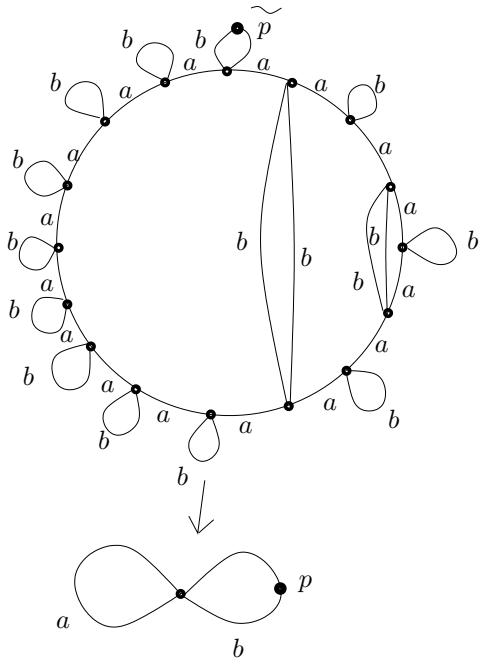}}\hspace{10mm}}
\caption{A cover of $G$, with $2n = 16$ and $\sigma = (1,2)(3,4)$}
\label{fig:1}
\end{figure}

 Let $\rho:T \rightarrow \tilde{G}/\tilde{p}$ be the universal
 cover. We give $T$ the standard path metric (see above).
 Let $\mathcal{E}_1$ be the set of edges of $T$ whose endpoints
 are equidistant from $\partial T$.
 Note that $\mathcal{E}_1$ is preserved under any isomorphism of $T$,
 and that $\mathcal{E}_1$ contains $\rho^{-1}$ of any 1-cycle.
 Let $\mathcal{E}_2$ be the set of edges in
 $\mathcal{E}(G) - \mathcal{E}_1$
 which are adjacent to at least two elements of $\mathcal{E}_1$.
 Then $\mathcal{E}_2$ is also isomorphism-invariant.
 Since no edge labeled $a$ has endpoints which are equidistant
 from $\partial (\tilde{G}/\tilde{p})$, then
 $\rho^{-1} \pi^{-1} a \cap \mathcal{E}_1 = \emptyset$, and
 since every edge labeled $a$ in $\tilde{G}$ is adjacent
 to a 1-cycle, then $\rho^{-1} \pi^{-1} a \subset \mathcal{E}_2$.
 Also, any element of $\rho^{-1} \pi^{-1} b - \mathcal{E}_1$ can be adjacent to
 (and distinct from) at most one element of $\mathcal{E}_1$,
 so $\rho^{-1} \pi^{-1} b \cap \mathcal{E}_2 = \emptyset$.
 Therefore $\mathcal{E}_2 = \rho^{-1} \pi^{-1} a$.

 Let $d_a: \mathcal{V}(T) \rightarrow \mathbb{Z}$
 be the function defined by setting $d_a(x)$ to be the distance
 from $x$ to $\partial T$
 along paths whose edges (except for the final edge) are all labeled $a$.
 Let $\mathcal{V}_i$ be the set of vertices of $T$
 with $d_a(x) = i$. Since $\mathcal{E}_2=\rho^{-1}\pi^{-1} a$
 is isomorphism-invariant,
 then $\mathcal{V}_i$ is isomorphism invariant.
 The transposition $(i j)$ occurs in the cycle decomposition
 of $\sigma$ if and only if there is an edge between
 $\mathcal{V}_{2 i}$ and $\mathcal{V}_{2 j}$, and
 thus the permutation $\sigma$ is determined by the isomorphism
 type of $T$.  Since the number of such permutations grows
 as a factorial in $n$, the number of inequivalent covers
 $(\tilde{G},\tilde{p})$ also grows as a factorial in $n$.
\end{proof}

\begin{lem} \label{initcover}
 Let $S$ be a hyperbolic surface, let $L$ be a
 non-separating collection of disjoint loops in $S$, and let $u$ be an
 arbitrary function of two variables.
 Then there is a finite cover $\tilde{S}$ of $S$,
 and a bouquet of two circles $G$ such that
\begin{enumerate}
\item there is a $(k,c)$-quasi-isometry
 $g:\tilde{S} \rightarrow G$,
\item both edges of $G$ have length greater
 than $u(k,c)$, and
\item $g(\tilde{L})$ is a point, for some 1-1 lift $\tilde{L}$ of $L$.
\end{enumerate}
\end{lem}

\begin{figure}[ht!]\footnotesize
\psfraga <-1pt, 0pt> {a}{$a$}
\psfraga <-1pt, -1pt> {b}{$b$}
\psfraga <0pt, 0pt> {r}{$\rho$}
\psfraga <0pt, 0pt> {S}{$S$}
\psfraga <-2pt, 0pt> {L}{$L$}
\psfraga <-5pt, 1pt> {fL}{$f(L)$}
\epsfxsize=.99\hsize \centerline{\epsfbox{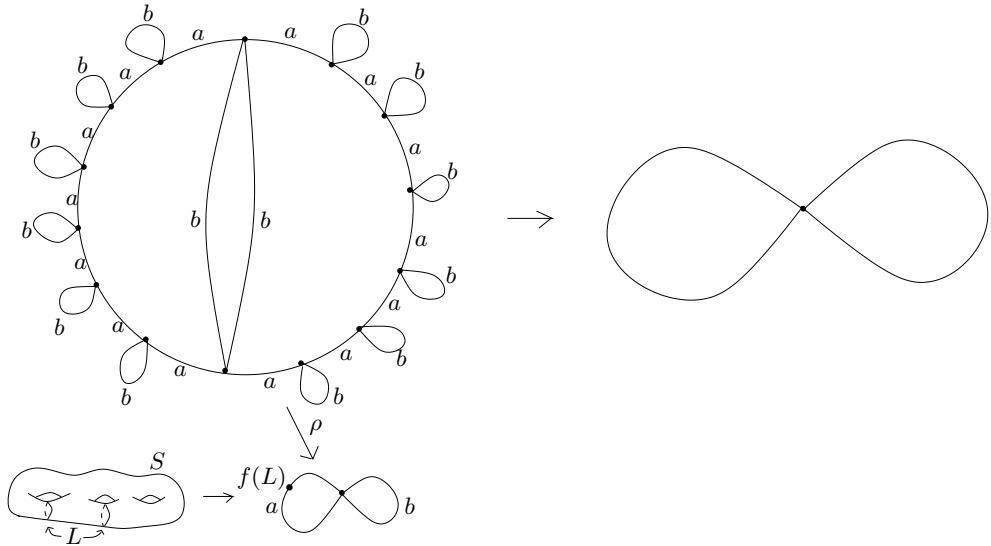}}
\caption{There is a cover of S which
 maps onto a bouquet of circles, with long edges.}
\label{fig:2}
\end{figure}

\begin{proof}
 After passing to a preliminary cover, we may assume that
 $S-L$ has positive genus. Then there is a
 map of $S$ onto a bouquet of two circles,
 $g_1: S \rightarrow G_1 = a \cup b$,
 such that $g_1$ induces a surjection
 on fundamental groups, and $g_1(L)$ is a point.
 We give $G_1$ the standard metric (so $a$
 and $b$ have length 1), and 
 since $S$ and $G_1$ are compact, $g_1$
 is a $(k,c_1)$-quasi-isometry for some $k,c_1$.

 Let $c=c_1+2$, and let $s > u(k,c)$.  Let
 $\pi:\tilde{G_1} \rightarrow G_1$
 be a $2s$-fold cover of $G_1$, where $\rho^{-1} a$
 is a single cycle, and $\pi^{-1} b$
 is a set of $2s-2$ cycles, as indicated in Figure \ref{fig:2}.
 There is an induced cover $\tilde{S}$ of $S$, to which $L$ lifts.
 Let  $\tilde{b}$ be the component of $\pi^{-1} b$
 which maps 2 to 1 onto $b$, and let $G$ be
 the graph $\pi^{-1} a$, with
 the vertices of $\tilde{b}$ identified to a point.
 Then $G$ is homeomorphic to a bouquet of two circles, and
 there is a $(1,2)$-quasi-isometry
 $g_2:\tilde{G_1} \rightarrow G$.
 Thus $g_2g_1: \tilde{S} \rightarrow G$ is a $(k,c)$-quasi-isometry,
 the edges of $G$ all have length at least $s > u(k,c)$,
 and there are 1-1 lifts $\tilde{L}$ of $L$ such that
 $g_2g_1 \tilde{L}$ is a point.\end{proof}

 Let $S$ be a closed hyperbolic surface, with a
 non-separating collection of loops $L$ in $S$, and let
 $(\tilde{S}, \tilde{L})$
 be a finite cover of $S$, with a 1-1 lift $\tilde{L}$ of $L$.
 Let $\tilde{S}/\tilde{L}$
 be the metric completion of $\tilde{S}-\tilde{L}$, which
 is a hyperbolic surface with geodesic boundary,
 and let $P \subset \mathbb{H}^2$ be the universal cover
 of $\tilde{S}/\tilde{L}$.
 For another such cover $(\tilde{S}^{\prime}, \tilde{L}^{\prime})$,
 we define analogous objects $\pi^{\prime}, \tilde{S}^{\prime}$
 and $P^{\prime}$.
 We say that the pairs $(\tilde{S}, \tilde{L})$ and
 $(\tilde{S}^{\prime}, \tilde{L}^{\prime})$ are \textit{inequivalent}
 if $P$ is not isometric to $P^{\prime}$.
\eject

\begin{lem} \label{inequiv}
 Let $(S,L)$ be as above.  Then the number of inequivalent covers
 $(\tilde{S}, \tilde{L})$ of degree $n$ is at least
 $(cn)!$ for some constant $c>0$.
\end{lem}

\begin{proof}
 We replace $S$ with the cover given by Lemma \ref{initcover}.
 We are free to do this, since the effect on
 covering degrees is linear.  We have a
 $(k,c)$-quasi-isometry $g:S \rightarrow G$,
 where $G$ is a bouquet of two circles,
 whose lengths are at least $u(k,c)$, where
 $u$ is a function of $k$ and $c$ which we will determine later.

 Let $(\tilde{G}, \tilde{p})$,
 and $(\tilde{G}^{\prime}, \tilde{p}^{\prime})$ be
 inequivalent covers of $G$, with
 induced covers $(\tilde{S}, \tilde{L})$
 and $(\tilde{S}^{\prime}, \tilde{L}^{\prime})$ of $S$.
 Let $P$ (resp. $P^{\prime}$)
 be the universal cover of $\tilde{S}/\tilde{L}$
 (resp. $\tilde{S}^{\prime}/\tilde{L}^{\prime}$),
 and let $T$ (resp. $T^{\prime}$)
 be the universal cover of $\tilde{G}/\tilde{p}$
 (resp. $\tilde{G}^{\prime}/\tilde{p}^{\prime}$).
 There is a $(k,c)$-quasi-isometry $h:P \rightarrow T$
 (resp. $h^{\prime}: P^{\prime} \rightarrow T^{\prime}$)
 obtained by lifting $g$ to $\mathbb{H}^2$ and then restricting the domain.
 Let $h^{\prime *}:T^{\prime} \rightarrow P^{\prime}$
 be a quasi-inverse of $h^{\prime}$.
 If $(\tilde{S},\tilde{L})$ and $(\tilde{S}^{\prime}, \tilde{L}^{\prime})$
 are equivalent, then there is an isometry $i:P^{\prime} \rightarrow P$.
 Then the map $hih^{\prime *}:T^{\prime} \rightarrow T$ is a
 $(k^{\prime},c^{\prime})$-quasi-isometry
 where $k^{\prime}$ and $c^{\prime}$ depend only on $k$ and $c$,
 and not on our choices of covers. By Lemma \ref{initcover},
 we may assume that the edges of $G$,
 and hence also of $T$ and $T^{\prime}$, 
 are longer than the constant $u(k^{\prime},c^{\prime})$
 given by Lemma \ref{quasi}.
 Then $T$ and $T^{\prime}$ are isomorphic, contradicting the assumption
 that $(\tilde{G}, \tilde{p})$ and
 $(\tilde{G}^{\prime}, \tilde{p}^{\prime})$ are inequivalent.
 So $(\widetilde{S}, \widetilde{L})$ and $(\widetilde{S}^{\prime},
 \widetilde{L}^{\prime})$ are inequivalent, and the result now
 follows from Lemma \ref{graphinequiv}.
\end{proof}

\section{Constructing immersed surfaces by gluing convex manifolds}
 The material in this section is based on Thurston's bending deformation,
 and the treatment of this which is given in \cite{BM}.
 We begin with a lemma about gluing convex manifolds.
 For a space $X$, we let $\mathring{X}$ denote the interior
 of $X$.

\begin{lem} \label{gluing}
Suppose that $X$ and $X^{\prime}$ are hyperbolic $n$-manifolds with
 convex boundary,
 and that $U, U^{\prime}$ are compact $n$-submanifolds,
 with an isometry $f:U \rightarrow U^{\prime}$
 such that
\begin{enumerate}
\item $f(\partial U \cap \mathring{X}) \subset \partial X^{\prime}$,
  $f^{-1}(\partial U^{\prime} \cap \mathring{X^{\prime}}) \subset
  \partial X$, and
\item $\partial X \cap f^{-1} \partial X^{\prime}$
  is a non-empty subsurface of $\partial X$.\\ Then the identification
  space $Y = X \cup_{U=fU} X^{\prime}$ is a hyperbolic manifold with
  convex boundary.
\end{enumerate}
\end{lem}

\begin{proof}
Condition 1 guarantees that the identification space $Y$ is a manifold,
 with $\partial Y = (\partial X - U) \cup (\partial X^{\prime} - U^{\prime})
 \cup (\partial X \cap f^{-1} \partial X^{\prime})$.
 Condition 2 guarantees that $\partial X - U$
 (resp. $\partial X^{\prime} - U^{\prime}$) can be extended to an open set 
 $V \subset \partial Y$ (resp. $V^{\prime} \subset \partial Y$),
 such that $Y$ is locally convex at all points in $V$.
 Then $V, V^{\prime}$ and
 $interior(\partial X \cap f^{-1} \partial X^{\prime})$
 define an open covering of $\partial Y$ by convex sets,
 and so $Y$ is locally convex at all points in $\partial Y$.
\end{proof}

 If $W$ is a hyperbolic 3-manifold with convex boundary, let $\Core(W)$
 denote the convex core of $W$.
 This is the unique smallest convex sub-manifold of $W$ which
 carries $\pi_1 W$.
 For any set $U$ in $W$, let
 $N_{\epsilon}(U) = \{ x \in U | d(x,U) \leq \epsilon \}$.

 Suppose $S$ is an orientable (possibly disconnected) hyperbolic
 surface with geodesic boundary,
 that $\ell_{ij} \subset \partial S$ ($i=1,2$, $j=1, ..., n$)
 are distinct boundary components, with orientations induced from $S$,
 and that $\length(\ell_{1j}) = \length(\ell_{2j})$ for all $j$. Let
 $L$ be the disjoint union of the $\ell_{ij}$'s, and let
 $S(L)$ be the surface formed by identifying each
 $\ell_{1j}$ and $\ell_{2j}$
 (preserving orientations) to a single loop, which we call $\ell_j$;
 we view $S(L)$ as a metric space, with
 path metric inherited from the $S_i$'s.
 For each $j$, let $\theta_j$ be angle between 0 and $\pi$.
 
\begin{lem} \label{bm}
 There is a number $r = r(\theta_1, ..., \theta_n) > 0$ such that, if 
 $L$ has a collar neighborhood of radius $r$ in $S$,
 then there is a hyperbolic 3-manifold $W$, with convex boundary,
 such that:
\begin{itemize}
\item[{\rm(a)}] $W = \Core(W)$,
\item[{\rm(b)}] $W \cong S(L) \times I$, and there is an
 isometric embedding $i: S(L) \rightarrow W$,
 such that $W$ retracts onto $i(S(L))$,
\item[{\rm(c)}] there is a regular neighborhood  $N(\ell_j)$ of each
 $\ell_i$ in $S(L)$ so that $i(N(\ell_j))$ is the union of
 a pair of geodesic annuli meeting at an angle $\theta_j$,
\item[{\rm(d)}] for some $\epsilon > 0$, the set
 $\mathcal{A} = \{ i(\ell_1), ..., i(\ell_n) \}$
 has the following geometric characterization:\newline
 $\mathcal{A} = \{$ infinite or closed geodesics $\gamma$ in $W|
 d(\gamma, F) > \epsilon$ for some component $F$ of $\partial W \}$.
\end{itemize}
\end{lem}

\begin{proof}
 We will assume that $n=1$, the proof in the general case being similar.
 To reduce notation, let $\ell_i = \ell_{i1}$, $\theta = \theta_1$,
 and $r=r_1$.

 The outline of the construction is as follows.
 We first construct a hyperbolic solid torus $V$,
 which is a convex thickening of two annuli.
 We then construct a convex hyperbolic 3-manifold $Z$,
 by gluing together $S(L) \times I - N(\ell_1)$,
 $S(L) \times I - N(\ell_2)$ and $V$.
 Then finally we let $W = \Core(Z)$.

 To construct the hyperbolic solid torus $V$, let $\epsilon > 0$,
 let $\tilde{\ell}$ be a geodesic in $\mathbb{H}^3$,
 let $Q_1$ and $Q_2$ be geodesic half-planes in $\mathbb{H}^3$,
 with $\partial Q_i = \tilde{\ell}$, and let
 the angle of intersection of $Q_1$ and $Q_2$ be $\theta$.
 Let $Q = Q_1 \cup Q_2$, let $B = N_{\epsilon}(Q)$, and let
 $H_Q$ be the convex component of $\mathbb{H}^3 - Q$.
 For $p \in \partial N_{\epsilon}(Q_i)$, let $P_p$
 be a hyperbolic plane such that $P_p \cap N_{\epsilon}(Q_i) = p$.

\medskip
\textbf{Claim}\qua{\sl There is an $r = r(\theta)$ such that,
 if $p \in \partial N_{\epsilon}(Q_i)$
 and $d(p, \tilde{\ell}) > r$, then $P_p \cap B = \{p\}$.}

\proof[Proof of Claim]
 There are two components of
 $\partial N_{\epsilon}Q_i - N_{\epsilon}(\tilde{\ell})$;
 one is contained in $H_Q$, which we denote
 $\partial^- N_{\epsilon}(Q_i)$, and the other is denoted
 $\partial^+ N_{\epsilon}(Q_i)$. 

 Without loss of generality, let $p \in \partial N_{\epsilon}(Q_1)$.
 Suppose $p \in \partial^+ N_{\epsilon}(Q_1)$.
 Let $P_1$ be the geodesic plane containing $Q_1$,
 and let $H_1$ be the half-space bounded by $P_1$ and containing
 $p$. If $r > \epsilon$, then
 $P_p \cap N_{\epsilon} Q_1 = P_p \cap N_{\epsilon}(P_1) = \{ p \}$, and
 so $P_p \subset H_1$.
 Since $N_{\epsilon}(Q_2) \cap H_1 \subset N_{\epsilon}(\tilde{\ell})$,
 and since $d(p, \tilde{\ell}) > \epsilon$,
 then $P_p \cap N_{\epsilon}(Q_2) = \emptyset$.
 Therefore $P_p \cap B = \{ p \}$.

 Suppose $p \in \partial^- N_{\epsilon}(Q_1)$.
 Let $Q^{\perp}$ be a hyperbolic plane which is orthogonal
 to $\tilde{\ell}$, and which contains $p$.
 Let $\alpha_i = Q_i \cap Q^{\perp}$, and let $\beta$
 be a geodesic ray bisecting the $\alpha_i$'s.
 Let $\delta(t)$ be a geodesic intersecting $\beta$
 orthogonally at a distance $t$ from the endpoint of $\beta$
 (see Figure \ref{fig:3}).

 Let $\partial^{\pm} N_{\epsilon}(\alpha_i)
            = \partial^{\pm} N_{\epsilon} Q_i \cap Q^{\perp}$.
 If $\delta(t) \cap B \neq \emptyset$, then
 by symmetry, $\delta(t)$ makes the same angle with both
 $\partial^- N_{\epsilon}(\alpha_1)$ and $\partial^- N_{\epsilon}(\alpha_2)$;
 call this angle $\psi(t)$.
 For small $t$, we have $\psi(t) > 0$, and for large enough $t$,
 we have $\delta(t) \cap N_{\epsilon}(\alpha_i) = \emptyset$,
 so there must be some $t_0$ for which $\psi(t_0) = 0$,
 and there is a corresponding geodesic $\delta(t_0)$ 
 which is tangent to $B \cap Q^{\perp}$ at two points.
  There is then a hyperbolic plane, $Q(t_0)$,
 orthogonal to $Q^{\perp}$, which contains $\delta(t_0)$,
 and which is tangent to $B$ at two points.

\begin{figure}[ht!]\small
\psfraga <-2pt, 0pt> {b}{$\beta$}
\psfraga <-5pt, 0pt> {d}{$\delta(t_0)$}
\psfraga <-3pt, 0pt> {a1}{$\alpha_1$}
\psfraga <-4pt, 0pt> {a2}{$\alpha_2$}
{\epsfxsize=2.5in \centerline{\epsfbox{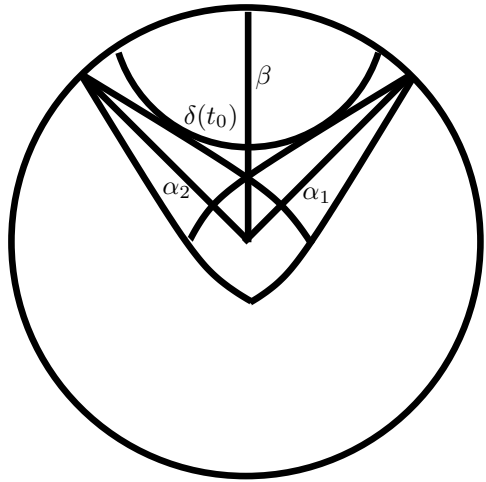}}\hspace{10mm}}
\caption{A convex thickening}\label{fig:3}
\end{figure}
 
 Let $r = d(\tilde{\ell}, \delta(t_0) \cap B)$, and suppose
 $d(p, \tilde{\ell}) > r$.
 Then $\delta(t_0)$ separates $p$ from $N_{\epsilon}(\alpha_2)$
 in $Q^{\perp} - B$.  Since the geodesic
 $P_p \cap Q^{\perp}$ can intersect $\delta(t_0)$ at most once,
 we must have $P_p \cap Q^{\perp} \cap N_{\epsilon}(\alpha_2) = \emptyset$,
 so $P_p \cap Q^{\perp} \cap B = \{ p \}$.
 Let $H^{\pm}$ be the half-spaces in $\mathbb{H}^3$ bounded by $Q^{\perp}$.
 By symmetry, if $P_p \cap N_{\epsilon}(Q_2) \neq \emptyset$, then
 $P_p \cap H^+ \cap N_{\epsilon}(Q_2)$ and
 $P_p \cap H^- \cap N_{\epsilon}(Q_2)$ are both non-empty, and
 therefore, by convexity, $P_p \cap Q^{\perp} \cap N_{\epsilon}(Q_2)
 \neq \emptyset$, which is a contradiction.
 Thus we conclude that $P_p \cap N_{\epsilon}(Q_2) = \emptyset$,
 so $P_p \cap B = \{ p \}$, proving the claim.\endproof

 Let $\Hull(B)$ denote the convex hull of $B$ in $\mathbb{H}^3$,
 which is the intersection of all closed half-spaces containing $B$.
 Let $\tau$ be a hyperbolic isometry (i.e. with real trace)
 in $\Isom^+(\mathbb{H}^3)$
 such that $\tilde{\ell}/\tau$ is isometric to $\ell$.
 Since $\tau$ is hyperbolic,
 it preserves both $Q_1$ and $Q_2$, and hence also $\Hull(B)$.
 Let $s > r(\theta)$, and let $V = [\Hull(B) \cap N_s(\tilde{\ell})]/\tau$.
  
 Let $S \times [-\epsilon, \epsilon]$ be a radius $\epsilon$ thickening
 of $S$ (more precisely, embed $\pi_1 S$ as a Fuchsian subgroup
 of $\Isom(\mathbb{H}^3)$, embed $\tilde{S}$ in $\mathbb{H}^3$ isometrically
 as a $\pi_1 S_i$-equivariant subset, and then define
 $S \times [-\epsilon, \epsilon] = N_{\epsilon}(\tilde{S})/\pi_1 S$).
 Let $U_i  = (V - N_r(\tilde{\ell})) \cap (N_{\epsilon}(Q_i)/\tau)$,
 and let
 $U_i^{\prime}=S \times [-\epsilon, \epsilon] \cap (N_s(\ell_i)- N_r(\ell_i))$.

 By the claim, $\Hull(B) - N_r(\tilde{\ell}) = B - N_r(\tilde{\ell})$,
 which implies that $U_i$ is isometric to $U_i^{\prime}$.
 By Lemma \ref{gluing}, we may form a hyperbolic manifold, $Z$, 
 with convex boundary, from
 $(S \times [-\epsilon, \epsilon]) - (N_r(\ell_1) \cup N_r(\ell_2)) \amalg V$
 by gluing each $U_i$ to $U_i^{\prime}$ via isometry.

 Since $Z$ has convex boundary,
 there is a unique smallest convex submanifold,
 $W = \Core(Z) \subset Z$, which is homotopy equivalent to $Z$.
 By the construction of $Z$, we have an isometric embedding
 $i:S(L) \rightarrow Z$, and $i(\ell)$ has a neighborhood in $Z$
 which is a union of geodesic annuli meeting at an angle $\theta$.
  Since $i(S)$ is totally geodesic in $Z$,
 with geodesic boundary, we have that $i(S) \subset W$,
 and so $i(S(L)) \subset W$.
 
 Since $\Core(\pi_1 W)$ is compact, then $\pi_1 W$ is a
 quasi-Fuchsian surface group, and
 the only possibilities are that $W \cong S(L) \times I$
 or that $W \cong S(L)$.
   Since $i(S(L))$ is not totally geodesic, $W \not\cong S(L)$,
 so $W \cong S(L) \times I$.  Thus $W$ satisfies properties $a,b$ and $c$
 of the lemma.

 We now prove that $W$ satisfies property $d$.
 The surface $i(S(L))$ is an embedded copy of $S(L)$ in $W$,
 which separates $W$ into two components, $W^+$ and $W^-$.
 Let $W^-$ be the component which is convex in a neighborhood
 of $\ell$, and let $F = W^- \cap \partial W$.

 The distance from $\ell$ to $F$ may be approximated to
 arbitrary precision by lengths of geodesic segments
 $\beta$ which are perpendicular to $\ell$.
 Let $\beta$ be such a segment, and let $\tilde{\beta}$
 and $\tilde{\ell}$ be lifts to $\tilde{W} \subset \mathbb{H}^3$.
 Let $A_1$,$A_2$ be radius $r$, totally geodesic annuli
 in $S(L)$, with $A_1 \cap A_2 = \ell$,
 and let $\tilde{A}_i \subset \tilde{W}$ be a lift of $i(A_i)$
 so that $\tilde{A}_1 \cap \tilde{A}_2 = \tilde{\ell}$.
 Let $Q^{\perp}$ be a geodesic plane containing $\tilde{\beta}$
 and perpendicular to $\tilde{\ell}$,
 and let $q = Q^{\perp} \cap \tilde{\ell}$.
 Then $Q^{\perp} \cap \tilde{A}_i = \alpha_i$ is a geodesic segment
 of length at least $r$, and $\alpha_1$ and $\alpha_2$
 meet at the point $q$ in an angle $\theta$.
 Let $T$ be the triangle with sides $\alpha_1$ and $\alpha_2$,
 and let $\delta$ be the third side of $T$.
 By convexity, $T \subset \tilde{W}$, and so
 $\length(\beta) \geq d(q, \delta)$.
 Let $T^{\prime}$ be a triangle with two ideal vertices
 and regular vertex $q$, with angle $\theta$, and
 let $s(\theta)$ be the distance from $q$ to the opposite side of $T^{\prime}$.
 Then as $r$ becomes large, $d(q, \delta)$ approaches $s(\theta)$.
 Therefore $d(\ell, F) \approx \length(\beta)  \approx  s(\theta)$,
 and these approximations can be made arbitrarily close
 by increasing $r$. Therefore, after possibly increasing $r$
 and decreasing $\epsilon$, we have that $d(\ell, F) > \epsilon$.

 Let $\gamma$ be some other geodesic in $W$.
 Then $\gamma$ also represents a geodesic in $Z$.
 Since $V$ contains a unique closed geodesic,
 and no infinite geodesics, then $\gamma \cap (Z - V) \neq \emptyset$.
 However, $(Z-V) \cong (S(L) \times [-\epsilon, \epsilon]) - N_r(\ell)$,
 and so $\gamma$ is within $\epsilon$ of both boundary
 components of $W$.
\end{proof}

 Next, some more notation. Suppose $f: S \rightarrow M$
 is a totally geodesic immersion,
 and suppose $\ell \subset f(S)$ is a double curve,
 along which the angle of self-intersection is $\theta$.
 Let $\ell_1, \ell_2 \subset S$ be two components of $f^{-1} \ell$,
 and let $L = \ell_1 \cup \ell_2$.

 Suppose that $\ell_i$ has a collar neighborhood of radius $r = r(\theta)$,
 where $r$ is the constant in Lemma \ref{bm}.
 Let $N(\ell_i)$ be an embedded collar neighborhood of $\ell_i$,
 and let $\partial N(\ell_i) = \{ \ell_{i1}, \ell_{i2} \}$.
 Further suppose, for convenience,
 that $L$ is non-separating. Let $S^-$ be a compact surface homeomorphic
 to $S - int N(L)$.  The hyperbolic structure
 on $S$ induces a hyperbolic structure on $S^-$ for which the components of
 $\partial S^- = \{ \ell_{11}, \ell_{12}, \ell_{21}, \ell_{22} \}$
 are all geodesics with the same length that $\ell$ has in $S$.
 Also, the orientation on $S_i$ induces an orientation on $\ell_{ij}$. 
 Let $S(L)$ be the surface obtained from $S^-$ by identifying $\ell_{1j}$ and
 $\ell_{2j}$, so that their orientations agree,
 and let $W = W(L, \theta) \cong S(L) \times I$
 be the hyperbolic 3-manifold given by Lemma \ref{bm}.
 The immersion $f:S \rightarrow M$ induces an
 immersion $f:S^- \rightarrow M$, which induces an
 immersion $g:S(L) \rightarrow M$, and extends to an immersion
 $g:W \rightarrow M$.

\begin{lem} \label{inject}
If $\ell_1$ and $\ell_2$ have collar
 neighborhoods of radius $r(\theta)$ in $S$, then the immersion
 $g:S(L) \rightarrow M$ is $\pi_1$-injective.
\end{lem}
 
\begin{proof}
 Since $W$ has convex boundary, the map $g:W \rightarrow M$
 is $\pi_1$-injective, so $g:S(L) \rightarrow M$ is
 $\pi_1$-injective also.
 \end{proof}

Lemma \ref{inject} provides a way of constructing a great
 many $\pi_1$-injective immersions.
 For suppose $\tilde{S}$ is a cover of $S$, and suppose 
 that $\ell$ has two (1-1) lifts, $\tilde{\ell}_1$ and $\tilde{\ell}_2$
 in $\tilde{S}$. We say that the pair
 $(\tilde{S}, \tilde{L}=\tilde{\ell}_1 \cup \tilde{\ell}_2)$
 is \textit{admissable} if the $\tilde{\ell}_i$'s
 have collar neighborhoods of radius $r$ in $\tilde{S}$.
 If $(\tilde{S}, \tilde{L})$ is admissable, then by Lemma \ref{inject},
 the immersion $g:\tilde{S}(\tilde{L}) \rightarrow M$ is $\pi_1$-injective.
 Our goal is to show that different choices
 of admissable pairs will often result in inequivalent
 immersions; counting the number of choices for such pairs
 will then prove Theorem \ref{main}.

 Suppose, then, that $(\tilde{S}, \tilde{L})$ and
 $(\tilde{S}^{\prime}, \tilde{L}^{\prime})$ are admissable pairs.
 Let $W$ be the convex thickening
 of $\tilde{S}(\tilde{L})$ provided by Lemma \ref{bm},
 with map $i:\tilde{S}(\tilde{L}) \rightarrow W$,
 and let $\ell = i(\tilde{L}) \subset W$.
 Let $P$ be the universal cover
 of $\widetilde{S}^-$. 

 We let $C$ be the universal cover of $W$; note that $C$ is isometric
 to the convex hull of the limit set of  $g_* \pi_1 \tilde{S}(\tilde{L})$.
 We also define analogous objects $W^{\prime}$,
 $i^{\prime}$, $\ell^{\prime}$, $\pi^{\prime}$, $\tilde{S}^{\prime -}$,
 $P^{\prime}$ and $C^{\prime}$ for $(\tilde{S}^{\prime}, \tilde{L}^{\prime})$.

 Recall from the previous section that the pairs $(\tilde{S}, \tilde{L})$
 and $(\tilde{S}^{\prime}, \tilde{L}^{\prime})$
 are \textit{inequivalent} if $P$ and $P^{\prime}$ are not isometric.

\begin{lem} \label{non-isom}
 If $(\tilde{S}, \tilde{L})$ and $(\tilde{S}^{\prime}, \tilde{L}^{\prime})$
 are admissable, inequivalent pairs,
 then $C$ and $C^{\prime}$ are not isometric.
\end{lem}

\begin{proof}
 Suppose $\phi: C \rightarrow C^{\prime}$ is an isometry.
 Identify $P$ with a component of the pre-image
 of $\widetilde{S}^-$ in $C$.
 By Property d of Lemma \ref{bm},
 $\phi (\pi^{-1} \ell) = \pi^{\prime -1} \ell^{\prime}$,
 so $\phi(P)$ is a totally geodesic surface in $C^{\prime}$,
 bounded by a subset of $\pi^{\prime -1} \ell^{\prime}$.
 We claim that $\phi(P)$ is isometric to $P^{\prime}$.

 To see this, let $A$ be a separating, geodesic annulus embedded in
 $W^{\prime}$, containing $\ell^{\prime}$.
 Since $\phi(P) \cap \partial C = \emptyset$,
 each component of $\pi^{\prime -1} A \cap \phi(P)$ is an infinite
 geodesic; since each component of $\pi^{\prime -1} A$
 contains a unique infinite geodesic, then each component of
 $\pi^{\prime -1} A \cap \phi(P)$ is a lift of $\ell^{\prime}$.
 Let $Q^{\prime}$ be the closure of a component of
 $\phi(P) - \pi^{\prime -1} \ell^{\prime}$, and let $\delta$ be a geodesic
 arc connecting two boundary components of $Q^{\prime}$.
 Since
  $\phi(P) \cap \pi^{\prime -1} \ell^{\prime} 
          = \phi(P) \cap \pi^{\prime -1} A$,
 then  $\pi^{\prime} \mathring{\delta} \subset W^{\prime} - A$, and
 since $ W^{\prime} - A$ retracts onto the totally geodesic surface
 $i(\tilde{S}^{\prime}-\ell^{\prime})$,
 and $\pi^{\prime} \delta$ is a  geodesic, then
 $\pi^{\prime} \delta \subset i(\tilde{S}^{\prime -})$.
 Therefore we may fix a copy of
 $P^{\prime}$ in $C^{\prime}$ which contains at least two boundary
 components of $Q^{\prime}$.  This implies
 that $P^{\prime}$ and $Q^{\prime}$ are contained in the same geodesic
 plane $H$ in $C^{\prime}$.
 Both $P^{\prime}$ and $Q^{\prime}$ are
 the closure of the unique component of $H - \pi^{\prime -1} \ell^{\prime}$
 which contains $\mathring{\delta}$, and
 therefore $P^{\prime} = Q^{\prime}$.

 Finally, if $R^{\prime}$ is any other component of
 $\phi(P) - \pi^{\prime -1} (\ell^{\prime})$, then
 $\pi^{\prime} P^{\prime} = \pi^{\prime} Q^{\prime}
 = \pi^{\prime} R^{\prime} = \tilde{S}^{\prime -}$,
 which forces $Q^{\prime}$ and $R^{\prime}$ to
 meet at an angle $\theta$,
 contradicting the fact that $\phi(P)$ is totally geodesic.
 Therefore $\phi(P) = Q^{\prime} = P^{\prime}$.
 \end{proof}

\begin{proof} (Of Theorem \ref{main})
 We have already proved the upper bound, and we now
 prove the lower bound.
 We are given a totally geodesic immersion of a hyperbolic
 surface $f: S \rightarrow M$, with a loop of
 transverse intersection $\ell$, at an angle $\theta$.
 Let $r=r(\theta)$ be the constant given by Lemma \ref{bm}.
 Let $\ell_1$ and $\ell_2$ be components of $f^{-1} \ell$,
 and let $L = \{\ell_1 \cup \ell_2\}$.
 Since $\pi_1 S$ is LERF (by \cite{Sc}), 
 there is a finite cover $\pi: \tilde{S} \rightarrow S$
 in which $L$ lifts to a non-separating
 curve with collar neighborhood of radius $r$.
 Replacing the immersion $f:S \rightarrow M$
 with the immersion $f \pi:\tilde{S} \rightarrow M$,
 we may assume that $L$ is non-separating,
 with a collar neighborhood of radius at least $r$.
 Then any further cover
 $(\tilde{S}, \tilde{L})$ of $(S, L)$ (where $\tilde{L}$
 is a 1-1 lift of $L$) is admissable.
 
 By Lemma \ref{inequiv}, the number of admissable, inequivalent covers
 $(\tilde{S}, \tilde{L})$ of $(S,L)$
 grows factorially in the covering degree of $\tilde{S}$, which
 is proportional to the genus of the resulting immersed surfaces.
 Since the pairs are admissable, the corresponding immersions
 are $\pi_1$-injective by Lemma \ref{inject}.
 Since the pairs are inequivalent, Lemma \ref{non-isom} implies that 
 the corresponding subgroups of $\pi_1 M$ are inequivalent
 under the relation generated by commensurability and conjugacy.
 Theorem \ref{main} follows.
\end{proof}

\rk{Remark} The proof actually shows that the surface subgroups
 are inequivalent under the relation generated by commensurability
 and conjugacy \textit{in ${PSL_2(C)}$.}

\section{Example: Reflection orbifolds}
 Now we consider the special case of
 right-angled reflection orbifolds.
 These always have transversely immersed, totally geodesic
 surfaces (corresponding to any pair of adjacent faces) and so
 Theorem \ref{main} applies.
 However, in this case there are combinatorial proofs
 which give more explicit information.

 We first give a combinatorial construction of 
 immersed incompressible surfaces
 in right-angled reflection orbifolds.  This is a well-known
 construction, and the treatment here is similar to
 that of \cite{BM}. We will need the following definitions and notation.

\begin{figure}[!ht]\small
\psfrag {a}{(a)}
\psfrag {b}{(b)}
\epsfxsize=3.2in \centerline{\epsfbox{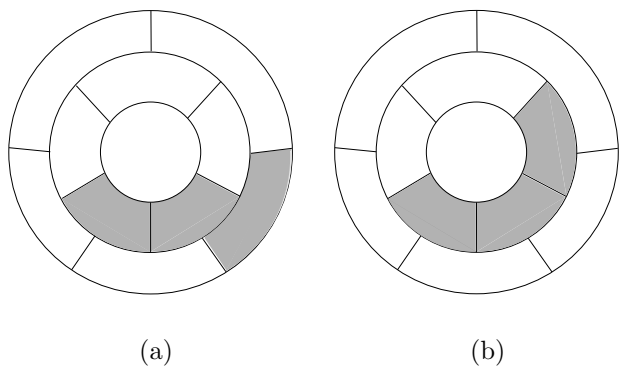}}
\caption{(a)\qua  A face disk satisfying the convexity condition\quad
 (b)\qua A face disk which does not satisfy the convexity condition}
\label{fig:4}
\end{figure}

Let $P$ be a right-angled polyhedron, let $\Gamma(P)$
 be the group generated by reflections in the faces of $P$,
 and let $\Gamma^+(P)$ denote the subgroup of orientation-preserving
 elements in $P$.  We define a \textit{face disk} of $P$ to be
 a connected and simply connected union of faces of $P$, and
 a \textit{face loop} on $P$
 to be a collection of $n$ faces $F_i$ such that
 $F_i$ is adjacent to $F_j$ if and only if $|i-j| = 1$ (mod $n$).
 Let $D$ be a face disk of $P$.
 We say that $D$ satifies the \textit{convexity condition} if
 the faces of $P$ transverse to $D$ form a face loop (see Figure \ref{fig:4}).

 If $F$ is a face in $P$, we let 
 $\rho_F :\mathbb{H}^3 \rightarrow \mathbb{H}^3$
 denote reflection in the plane containing $F$.
 If $X$ and $Y$ are spaces, and $\rho$ is a homeomorphism from a subset of $X$
 to a subset of $Y$, let $X \cup_{\rho} Y$ denote $X$ glued
 to $Y$ via the gluing map $\rho$.

 The incompressibility results of this section are essentially an elaboration
 of the following elementary fact: given a collection of planes in
 $\mathbb{H}^3$, each meeting the next at right angles,
 and such that non-consecutive planes are disjoint,
 then the group generated by their reflections is isomorphic to the
 obvious 2-dimensional reflection group.

\begin{lem} \label{convexlemma}
 Suppose that $P_1$ and $P_2$ are right-angled polyhedra in
 $\mathbb{H}^3$, that $D_i \subset P_i$ are face disks satisfying the
 convexity condition, and that $F_i \subset P_i$ are faces transverse to $D_i$,
 with $D_i \cap F_i = e_{i1}, ..., e_{im}$.  Suppose there is an
 orientation-reversing isometry $\rho: F_1 \rightarrow F_2$ such that
 $\rho(F_1, e_{11}, ..., e_{1m}) = (F_2, e_{21}, ..., e_{2m})$. 
 Then $D_1 \cup_{\rho} D_2$ satisfies the convexity condition
 in $P_1 \cup_{\rho} P_2$.
\end{lem}

\begin{proof}
 Since $P_i$ are right-angled polyhedra,
 $P_1 \cup_{\rho} P_2$ is a right-angled polyhedron also.
 
 The lemma then follows from the easily verified fact that
 the amalgamation of two face loops along a common
 face is again a face loop.
 \end{proof}

 Since $D$ abuts the reflection planes of $P$
 in right angles, it defines a closed 2-orbifold $\bar{D}$
 in $\mathbb{H}^3 / \Gamma(P)$ by making the edges
 of $D$ into reflectors. 
 The image of $\pi_1 \bar{D}$ in $\Gamma$
 is generated by the reflections in the faces of $P$ which
 $D$ interesects transversely.  The following result
 about the incompressibility of $\pi_1 \bar{D}$ is a reformulation
 of a well-known property of Coxeter groups.

\begin{pro} \label{incomp}
 If $D$ satisfies the convexity condition, then
 the fundamental group of $\bar{D}$ injects into $\Gamma(P)$.
\end{pro}

\begin{proof}
We fix a copy of $D$ and $P$ in $\mathbb{H}^3$.
 It will suffice to show that the union of $D$ with its translates under
 $\pi_1(\bar{D})$ forms an embedded disk.

 Since  $\bar{D}$ is a right-angled reflection orbifold, it may be developed
 in $\mathbb{H}^3$ by successively doubling $(P, D)$ across the faces
 of $P$ which $D$ intersects transversely.

 Let $F$ be a face of $P$ which $D$ intersects transversely,
 and let $P_2 = P \cup \rho_F P$.
 Then $D$ develops in $P_2$ as
 $D \cup \rho_F D$.  By Lemma \ref{convexlemma}, $D \cup \rho_F D$
 satisfies the convexity condition in $P \cup \rho_F P$.
 In particular, $D \cup \rho_F D$ is a disk.

 Iterating this process, we obtain after each doubling an embedded
 disk.  It follows that the union of $D$ with all its translates
 is embedded and simply connected.
\end{proof}

 Suppose that $\widehat{P} \subset \mathbb{H}^3$ is a right-angled
 polyhedron which is made up of copies of $P$,
 and that $D$ is a face disk in $\widehat{P}$
 satisfying the convexity condition.
 By Theorem \ref{incomp}, the reflection orbifold $\bar{D}$
 is incompressible in $\mathbb{H}^3/\Gamma(\widehat{P})$, so
 under the covering map
 $\pi: \mathbb{H}^3/\Gamma(\widehat{P}) \rightarrow \mathbb{H}^3/\Gamma(P)$,
 it projects to an immersed incompressible 2-orbifold
 $\pi \bar{D}$.
 As we did in Section 3, we can cut and paste $\pi \bar{D}$,
 to create many different
 immersions of 2-orbifolds in $\mathbb{H}^3/\Gamma(P)$.
 
 Let $e_1, e_2 \subset \partial D$ be edges
 of $\widehat{P}$, let $F_i$ be the faces adjacent to
 $e_i$ in $D$, and let $F_i^{\prime}$ be the face adjacent to $e_i$
 in $\widehat{P} - D$. Suppose there is an orientation-preserving element
 $\gamma \in \Gamma(P)$ such that
 $\gamma (F_1, e_1) = (F_2, e_2)$.
 We ``cut'' $\bar{D}$ by first making $e_1, e_2$ into boundary
 edges (not reflector edges); call this orbifold $\bar{D}^{\prime}$.
 We then ``paste'' by forming the quotient orbifold
 $\bar{D}^{\prime}/\rho_{F_2^{\prime}} \gamma$.
 Since $\rho_{F_2} \gamma$ takes $e_1$ to $e_2$, this is a closed
 2-orbifold, and since $\rho_{F_2} \gamma \in \Gamma(P)$,
 there is a natural map from $\bar{D}^{\prime}/\rho_{F_2} \gamma$
 into $\mathbb{H}^3/\Gamma(P)$.

 We repeat this procedure at different edges of $D$,
 and encode the collection of gluings by an order 2 permutation
 $\sigma$, defined on the edges of $D$.
 Let $D_{\sigma}$
 denote the immersed, closed 2-orbifold in $\mathbb{H}^3/ \Gamma(P)$
 obtained by cutting and pasting $\pi \bar{D}$ according to the
 permutation $\sigma$, and let
 $f_{\sigma}:D_{\sigma} \rightarrow \mathbb{H}^3/\Gamma(P)$ be the
 associated immersion.
 We define the \textit{orbifold degree} of a vertex $v$
 of $D_{\sigma}$ to be the degree of $v$ in the universal
 cover of $D_{\sigma}$.

  \begin{pro} \label{paste}
 If $D$ satisfies the convexity condition, and
 if each vertex of $D_{\sigma}$ has (orbifold) degree 4,
 then the immersion $f_{\sigma}:D_{\sigma} \rightarrow \mathbb{H}^3/\Gamma(P)$
 is $\pi_1$-injective.
\end{pro}

\begin{proof}
The proof is similar to the proof of Theorem \ref{incomp}.
 We fix a copy of $\widehat{P}$ and $D$ in $\mathbb{H}^3$
 and show that the translates of $D$ under $\pi_1(D_{\sigma})$
 form an embedded disk.

The translates of $D$ under $\pi_1(D_{\sigma})$
 may be obtained by successively doubling $\widehat{P}$
 via the various symmetries in $\Gamma(P)$ which correspond
 to the permutation $\sigma$.
 After the first doubling, then, $D$ is glued to $\rho D$
 along some edge, via some $\rho \in \Gamma(P)$.
 The faces $D$ and $\rho D$ satisfy the convexity condition,
 in $\widehat{P}$ and $\rho \widehat{P}$ respectively,
 so by Lemma \ref{convexlemma}, their union satisfies the convexity
 condition in $\widehat{P} \cup \rho \widehat{P}$.

 Doubling again along a face adjacent to the initial doubling face,
 we get four copies of $D$ around a vertex,
 which by the vertex condition will glue appropriately.
 Again, by Lemma \ref{convexlemma}, this union will be a face
 disk satisfying the convexity condition.

 Iterating this process, we obtain after each doubling an embedded disk, and
 therefore the union of $D$ with all its translates is an embedded disk.
\end{proof}

We now give an explicit construction of a family of inequivalent
 immersed surfaces in a right-angled reflection orbifold.

We are given a right-angled polyhedron $P$.
 Compact Coxeter polyhedra in $\mathbb{H}^3$
 have vertices of degree three, and no faces
 with fewer than five edges;
 an Euler characteristic computation then
 shows that every Coxeter polyhedron has a face with 
 exactly five edges.  

\begin{figure}[ht!]\footnotesize
\psfraga <-1pt,-1pt> {F1}{$F_1$}
\psfraga <-3pt,-1pt> {F2}{$F_2$}
\begin{center}
\epsfxsize 3.5in\epsfbox{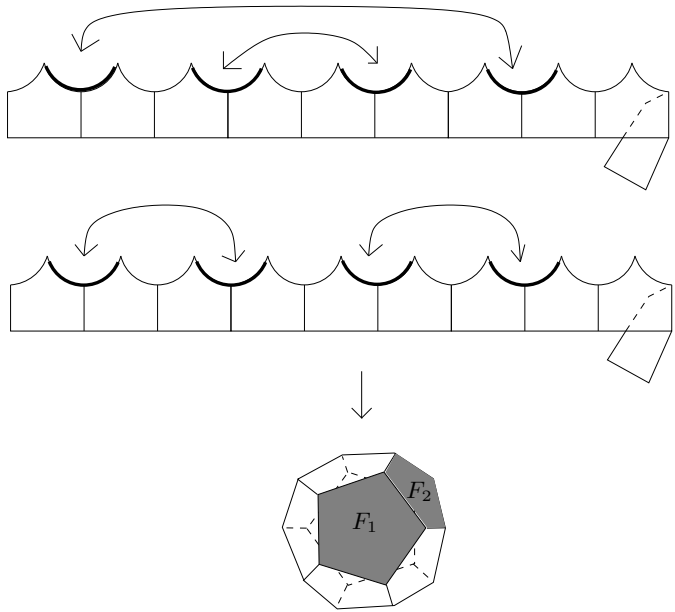}
\caption{Different gluings will result in inequivalent surfaces} 
\label{fig:5}\end{center}
\end{figure}

 Let $F_1$ be a face of $P$ with five edges. 
 We paste together $n$ copies of $F_1$, and attach one additional
 adjacent face $F_2$, as indicated in Figure \ref{fig:5};
 the resulting complex, which we call $D$, is a face disk
 for a right-angled Coxeter polyhedron $\widehat{P}$
 made up of $n$ copies of $P$.
 We cut and paste the reflection orbifold $\bar{D}$ along 
 alternating edges of $D$, as shown in Figure \ref{fig:5}.
  Associated to any order 2 permutation, $\sigma$, of these edges
 there is an edge pairing, and a resulting 2-orbifold $D_{\sigma}$,
 with an immersion $f_{\sigma}:D_{\sigma} \rightarrow \mathbb{H}^3/\Gamma(P)$.

\begin{lem} \label{incompex}
The immersion $f_{\sigma}:D_{\sigma} \rightarrow \mathbb{H}^3/\Gamma(P)$
 is $\pi_1$-injective.
\end{lem}

\begin{proof}
 The union of two adjacent faces of a compact Coxeter polyhedron
 always satisfies the convexity condition, since such a polyhedron
 cannot have face loops of length less than five.
 Therefore $D$ satisfies the convexity condition.
 Also, since $\sigma$ has order two,
 the vertices of $D_{\sigma}$ all have orbifold degree four,
 and so the lemma follows from Proposition \ref{paste}.
\end{proof}

The resulting immersions are in fact inequivalent.  The key point is that
 we have control over the limit sets.

Let $\tilde{P}$ be a fixed lift of $P$ to $\mathbb{H}^3$,
 and let $\tilde{D} \subset \tilde{P}$ be a fixed lift of $D$,
 which determines a fixed map
 $f_{\sigma *} \pi_1 D_{\sigma} \rightarrow \Gamma(P)$.
 Let $\Lambda(f_{\sigma *} \pi_1 D_{\sigma})$
 denote the limit set of $f_{\sigma *} \pi_1 D_{\sigma}$
 in $\mathbb{H}^3$.
 Let $\mathcal{H} = \{H_1, ... , H_n\} \subset \mathbb{H}^3$
 denote the half-spaces
 which intersect $\tilde{P}$ in faces orthogonal to $\tilde{D}$;
 we label these so that $H_i \cap H_{i+1} \neq \emptyset$,
 and $H_i \cap H_j = \emptyset$ if $|i-j| \geq 2$ (mod $n$).
Then we have:

\begin{lem} \label{limitset}$\phantom{99}$
\begin{enumerate}
\item $\Lambda(f_{\sigma*} \pi_1 D_{\sigma}) \cap H_i \neq \emptyset$
 for all $i$.\\
\item $\Lambda(f_{\sigma *} \pi_1 D_{\sigma}) \subset \bigcup_{i=1}^n H_i$.\\
\end{enumerate}
\end{lem}

\begin{proof}
 Each $H_i$ is stabilized by $<\rho_{H_{i-1}}, \rho_{H_{i+1}}>$.
 Since $H_{i-1}$ and $H_{i+1}$ are non-adjacent, this group is
 infinite, and so it has a limit point in $H_i$, proving 1.

 Since every point in $\Lambda(f_{\sigma *} \pi_1 D_{\sigma})$
 is a limit of a sequence of points in $\tilde{D}_{\sigma}$, then to prove 2,
 it is enough to show that 
 $\tilde{D}_{\sigma} \subset \bigcup_{i=1}^n H_i \cup \tilde{D}$.

 As in the proof of Lemma \ref{paste}, $\tilde{D}_{\sigma}$
 may be constructed by successive doublings.
 More precisely, there is a sequence of Coxeter polyhedra
 $\tilde{P} = P_1, P_2, ...$,
 and face disks $D_i \subset P_i$, such that
 $(P_{i+1}, D_{i+1})
 = (P_i, D_i) \cup (\gamma_i P_i, \gamma_i D_i)$,
 for some $\gamma_i \in f_{\sigma *} \pi_1 D_{\sigma}$,
 and  $\tilde{D}_{\sigma} = \bigcup_i D_i$.
 Let $\mathcal{Q}_i$ be the set of hyperbolic planes of $P_i$ which
 are transverse to $D_i$.

\medskip
\textbf{Claim}\qua For every $i$, 
 if $Q \in \mathcal{Q}_i$,then $Q \subset \bigcup_{H_j \in \mathcal{H}} H_j$.\\

\proof[Proof of claim]
 The statement is true for $i=1$, and suppose it is true
 for $i = n$. Let $J_n$ be the plane fixed
 by $\gamma_n$, and suppose without loss of generality that
 $J_n \subset H_1$.
 All of the planes in $\mathcal{Q}_{n+1} - \mathcal{Q}_n$
 correspond to faces of $\rho_n P_n$
 which are not adjacent to the face containing $J_n$.
 Since $\gamma_n P_n$ is a Coxeter polyhedron, these planes
 are disjoint from $J_n$,  and
 therefore $J_n$ separates every plane in $\mathcal{Q}_{n+1} - \mathcal{Q}_n$
 from $\tilde{D}$. Therefore every plane
 in  $\mathcal{Q}_{n+1} - \mathcal{Q}_n$ is contained in $H_1$,
 and so every plane in $\mathcal{Q}_{n+1}$ is contained in
 some $H_j$.  The claim follows by induction.
 
 From the claim it follows that $D_j \subset \bigcup_i H_i$
 for all $j$, and so $\tilde{D}_{\sigma} \subset \bigcup_i H_i$,
 as required.
\end{proof}

 Let $\sigma$ and $\sigma^{\prime}$ be distinct order two permutations.
 We choose arbitrary lifts $\tilde{D}$ and $\tilde{D}^{\prime}$
 of $D$ in $\mathbb{H}^3$, and these determine
 maps $f_{\sigma *}: \pi_1 D_{\sigma} \rightarrow \Gamma(P)$
 and $f_{\sigma^{\prime} *}: \pi_1 D_{\sigma^{\prime}} \rightarrow \Gamma(P)$.
 Let $\Gamma = f_{\sigma *} \pi_1 D_{\sigma}$, and
 $\Gamma^{\prime} = f_{\sigma \prime *} \pi_1 D_{\sigma^{\prime}}$.

\begin{pro}\label{refnonconj}
 No conjugate of $\Gamma$ is commensurable with $\Gamma^{\prime}$
 in $\Gamma(P)$.
\end{pro}

 Before proving this, we need one further lemma.
 Let $A, B \subset \mathbb{H}^3$, and 
 let $\mathcal{Q}$ be the collection of all the geodesic
 planes in $\mathbb{H}^3$ which project to $\partial P$.
 Define the \textit{combinatorial distance}
 $d_P(A, B)$ to be the
 number of planes in $\mathcal{Q}$ which separate $A$ and $B$.
 Let $\mathcal{H} = \{H_1, ..., H_k\}$
 be the half spaces perpendicular to $\partial \tilde{D}$. 
 For each $H_i$, there is an element $\gamma_i \in \Gamma$
 such that $\gamma_i \tilde{D} \cap \tilde{D} \subset \partial H_i$.
 
\begin{lem}\label{word}
Suppose $A \subset \mathbb{H}^3$, and that $H_i \supset A$.
 Then $d_P(\gamma_i \tilde{D}, A) < d_P(\tilde{D},A)$.
\end{lem}

\begin{proof}
 Otherwise there is a plane $Q \in \mathcal{Q}$ which
 separates $\gamma_i \tilde{D}$ from $A$ but
 does not separate $\tilde{D}$ from $A$.
 Then $Q$ must separate a subset of $\tilde{D}$
 from $\gamma_i \tilde{D}$, and so
 $Q \cap \mathring{\tilde{D}} \neq \emptyset$.
 Also, $Q \neq \partial H_i$, and so there is a $j$
 such that $Q, \partial H_i$, and $\partial H_j$ all intersect,
 which forces $Q$ to be adjacent to $H_i$, contradicting
 the assumption that $Q \cap \mathring{\tilde{D}} \neq \emptyset$.
\end{proof}
 
\begin{proof}[Proof of Proposition \ref{refnonconj}]
 Since $\Gamma$ and $\Gamma^{\prime}$
 were chosen only up to an arbitrary conjugation (depending
 on the choice of lifts $\tilde{D}$ and $\tilde{D}^{\prime}$),
 it is enough to show that $\Gamma$ is not commensurable
 with $\Gamma^{\prime}$, and for this, it is enough to show
 that $\Lambda(\Gamma) \neq \Lambda(\Gamma^{\prime})$.

 Suppose that  $\Lambda(\Gamma) = \Lambda(\Gamma^{\prime})$,
 and suppose also that $\tilde{D} \neq \tilde{D}^{\prime}$.
 Let $\mathcal{H} = \{H_1, ..., H_k\}$
 and $\mathcal{H}^{\prime} = \{H_1^{\prime}, ..., H_k^{\prime}\}$
 be the half spaces perpendicular to $\partial \tilde{D}$
 and $\partial \tilde{D}^{\prime}$, let $\gamma_i \in \Gamma$
 such that $\gamma_i \tilde{D} \cap \tilde{D} \subset \partial H_i$,
 and define $\gamma_i^{\prime}$ similarly.
 Then Lemma \ref{limitset} implies
 that every $H_i^{\prime}$ is contained in some $H_j$, and
 in particular $\tilde{D}^{\prime} \subset H_j$
 for some $j$. By Lemma \ref{word},
 $d_P(\gamma_j \tilde{D},\tilde{D}^{\prime})
 < d_P(\tilde{D}, \tilde{D}^{\prime})$, and
 we may replace $\tilde{D}$ by $\gamma_j \tilde{D}$,
 without affecting $\Lambda(\Gamma)$.
 Eventually, we have $d_P(\tilde{D}, \tilde{D}^{\prime}) = 0$,
 so $\tilde{D}$ and $\tilde{D}^{\prime}$ intersect in some lift of $P$. 

 Recall $\tilde{D}$ and $\tilde{D}^{\prime}$ both consist
 of $n$ lifts of $F_1$ and a single lift of $F_2$.
 If $\tilde{D}$ and $\tilde{D}^{\prime}$
 overlap on some lift of $F_1$ and differ on the lift of $F_2$,
 then there is some $H_j$ (adjacent to the lift of $F_2$)
 which is disjoint from every $H_i^{\prime}$,
 contradicting Lemma \ref{limitset}.
 So $\tilde{D}$ and $\tilde{D}^{\prime}$ contain the same lift of $F_2$, and 
 it then follows that $\tilde{D}$ and $\tilde{D}^{\prime}$ co-incide.

 We now have $\tilde{D} = \tilde{D}^{\prime}$,
 but by assumption, the gluings $\sigma$ and $\sigma^{\prime}$ disagree
 on some edge, which lifts to some $H_i$.
 Then $\mathcal{F} = \tilde{D} \cup \gamma_i \tilde{D}$ and
 $\mathcal{F}^{\prime} = \tilde{D} \cup \gamma_i^{\prime} \tilde{D}$
 are face disks in a Coxeter polyhedron $\widehat{P}$
 (made of $2n$ copies of $\tilde{P}$)
 which satisfy the convexity condition.
 Let $\mathcal{H}$ and $\mathcal{H}^{\prime}$ be the sets of half-spaces
 transverse to $\mathcal{F}$ and $\mathcal{F}^{\prime}$.
 Since $\sigma$ and $\sigma^{\prime}$ disagree,
 there is an $H_i \in \mathcal{H}$ which is disjoint from
 every $H_j \in \mathcal{H}^{\prime}$.  
 Lemma \ref{limitset} then
 shows that $\Lambda(\Gamma) \neq \Lambda(\Gamma^{\prime})$,
 which is a contradiction.
\end{proof}

 The lower bound given in Theorem \ref{reforbs} can be computed 
 from  Proposition \ref{refnonconj} as follows:
 first we must modify the construction of $D_{\sigma}$
 to produce closed, non-singular surfaces, instead of 2-orbifolds.
 This can be done by starting with 4 copies of $F_1$
 and $4n$ copies of $F_2$, which can be glued into a closed surface.
 There $4((n-1)/2) = 2n - 2$ edges (corresponding to the dark
 edges of Figure \ref{fig:5}) along which we may cut and paste
 via an arbitrary transposition $\sigma$, to get an immersion of a
 closed surface $S_{\sigma}$.
 By the same method as the proof of Proposition \ref{refnonconj},
 it can be shown that different choices of $\sigma$ result in inequivalent,
 $\pi_1$-injective immersions. 

  The Euler characteristic of $S_{\sigma}$
 is $-n-c$, for some fixed constant $c$ depending on the number
 of sides of $F_2$,
 and there are at least $(n-1)!$ choices of transpositions
 $\sigma$ resulting in inequivalent immersions.
 We have $g = (n+c)/2 + 1$, so $n-1 = 2g - c - 3$, 
 and so $s_1(M,g) > (2g-c-3)! > e^{g \log g}$, for large $g$.

Finally, we give a combinatorial proof of the upper bound in this setting.
  The first step is to show that every incompressible immersed surface
 in a reflection orbifold can be homotoped to a standard form.
 
 We view $\mathbb{H}^3/\Gamma(P)$ as an orbifold, with an underlying
 polyhedral complex determined by $P$.

\begin{lem} \label{normal}
Let $S$ be a closed surface of positive genus,
 and let $\Gamma \subset \Gamma(P)$, with $\Gamma \cong \pi_1 S$.
 Then there is a polyhedral 2-complex $K$ on $S$, and
 a immersion $f:S \rightarrow \mathbb{H}^3/\Gamma(P)$,
 which is cellular with respect to $K$ and $P$, such that
 $f_*(\pi_1 S) = \Gamma$, and such that the universal
 cover of $S$ lifts to an embedded plane
 in $\mathbb{H}^3$.
\end{lem}

\begin{proof}
There is a unique cover $\pi:M \rightarrow \mathbb{H}^3/\Gamma(P)$
 such that $Image(\pi_*) = \Gamma$. The pre-image of $P$ in $M$
 defines a polyhedral complex $K(P)$ on $M$.
 By Scott's compact core theorem, $M$ has a compact core $M^-$,
 such that $\pi_1(M^-) \cong \Gamma$, and Waldhausen's theorems,
 $M^- \cong S \times I$.  Let $K(P)^-$ denote the polyhedra
 of $K(P)$ which intersect $M^-$.  
 There is an embedding $i: S \rightarrow \partial K(P^-)$,
 and the polyhedral structure of $\partial K(P^-)$ pulls back
 to a polyhedral structure $K$ on $S$.
 Then $\pi i:S \rightarrow \mathbb{H}^3/\Gamma(P)$ is the desired immersion. 
\end{proof}

 Let $\tilde{P}$ be a lift of $P$ to $\mathbb{H}^3$.
 A sub-complex $\mathcal{F}$ of the 2-skeleton of $K$
 is a \textit{face disk} if $f|_{\mathcal{F}}$
 lifts to a face disk of $\tilde{P}$.
 If $v$ is a vertex of $\partial \mathcal{F}$, $deg_K(v)$
 denotes the degree of $v$ in $K$,
 whereas $deg_{\mathcal{F}}(v)$ denotes the degree of
 $v$ in $\mathcal{F}$.
  Let $\mathcal{V}(\mathcal{F},i)
     = \{ v \in \partial \mathcal{F}| deg_{\mathcal{F}}(v) = i\}$. 

 We require the following general fact about compact Coxeter polyhedra
 in $\mathbb{H}^3$.

\begin{lem} \label{degree}
Let $\mathcal{F}$ be a face disk of $P$.  If
 $|\mathcal{V}(\mathcal{F},2)| \leq 4$,
 then $\mathcal{F}$ contains more than half of the faces of $P$.
\end{lem}

\begin{proof}
 Let $\mathcal{F}^{\prime} = P - \mathcal{F}$.
 If $|\mathcal{V}(\mathcal{F},2)|=0$, then $\mathcal{F}^{\prime}$
 is a single face.
 If $|\mathcal{V}(\mathcal{F},2)|=1$, then
 $|\mathcal{V}(\mathcal{F}^{\prime},3)|=1$, which
 implies the existence of an edge in $P$ which
 is adjacent to only one face, which is impossible.
 If $|\mathcal{V}(\mathcal{F},2)|=2$, then
 $\mathcal{F}^{\prime}$ consists of exactly two faces of $P$.
  If $|\mathcal{V}(\mathcal{F},2)|=3$, then
 $\mathcal{F}^{\prime}$ consists of
 three faces meeting at a common vertex.
 If $|\mathcal{V}(\mathcal{F},2)|=4$, then, adjacent to its boundary
 $\mathcal{F}^{\prime}$ contains faces $F_1, ..., F_4$,
 with $F_i \cap F_{i+1} \neq \emptyset$ (mod 4).
 Since $P$ is a right-angled, compact Coxeter polyhedron,
 it contains no face loops of length three or four,
 and no vertices of degree four, and so this is impossible.

 We conclude that $\mathcal{F}^{\prime}$ contains
 at most three faces.  Since a compact, right-angled Coxeter polyhedron
 in $\mathbb{H}^3$ must contain at least twelve faces, the result follows.
 \end{proof}

 A face disk $\mathcal{F}$ of $K$ is \textit{special} if
\begin{enumerate}
\item $deg_K(v) \geq 4$ for all $v \in \partial \mathcal{F}$, and\\
\item $|\mathcal{V}(\mathcal{F},2)| \geq 5$.
\end{enumerate}

\begin{lem} \label{special}
 There is a choice of $f$ and $K$ as in Lemma \ref{normal},
 such that $K$ has a partition into special face disks.
\end{lem}
 
We say that the resulting immersed surface $S$ is in \textit{standard form}.
 For a complex $K$, let $K^{(i)}$ denote the $i$-skeleton,
 and let $|K^{(i)}|$ denote the cardinality of $K^{(i)}$.

\begin{proof}
Let $f$ and $K$ be as given by Lemma \ref{normal}.
 If $F_1, F_2$ are in $K^{(2)}$, say that $F_1 \equiv F_2$
 if there is a vertex $v \in F_1 \cap F_2$ with $deg_K(v) = 3$.
 Let $\mathcal{F}_1, ..., \mathcal{F}_n$
 be the equivalence classes of $\equiv$.
 Every $\mathcal{F}_i$ lifts to a single copy of $P$ in $\mathbb{H}^3$,
 and so $\mathcal{F}_i$ is a face disk for all $i$.
 
 Suppose $|\mathcal{V}(\mathcal{F}_i,2)| \leq 4$ for some $i$.
 Then we enlarge $\mathcal{F}_i$ to a maximal face disk,
 by adding as many faces of $K^{(2)}$ as possible.
 Let $\mathcal{F}_i^{\prime}$ denote $\partial P - f(\mathcal{F}_i)$.
 Then $f|_{\mathcal{F}}$ can be homotoped rel. boundary,
 and the polyhedral structure on $S$ can be changed,
 so that $f(\mathcal{F}_i) = \mathcal{F}_i^{\prime}$, and $f$ remains
 a cellular immersion.  Also, when lifted to the universal cover,
 this homotopy is supported on a single copy of $P$, and one
 may check that the lift of $f$ remains an embedding.  
 We then obtain new equivalence classes,
 and repeat the process.  By Lemma \ref{degree}, every step
 reduces $|K^{(2)}|$, so the process must eventually terminate
 in a special complex.
\end{proof}

Suppose now that $f:K \rightarrow \mathbb{H}^3/\Gamma(P)$ is as
 given by Lemma \ref{special}.
 Let $K^{\prime}$ be the corresponding complex of special face disks.

\begin{lem}\label{deg}
 If $v$ is a vertex of $K$ with $deg_K(v) \leq 5$,
 then $deg_{K^{\prime}}(v) \geq 4$.
\end{lem}

\begin{proof}
 Suppose $deg_{K^{\prime}}(v) \leq 3$.
 Let $\mathcal{F}_i$, $i=1, ..., n \leq 3$ be the face disks 
 of $K^{\prime}$ which are adjacent to $v$.
 Lift $f$ to an embedding $\tilde{f}:\tilde{S} \rightarrow \mathbb{H}^3$,
 choose a lift $\tilde{v}$ of $v$, and
 let $\tilde{\mathcal{F}}_i$ be the lift of $\mathcal{F}_i$
 adjacent to $v$.  Then $\tilde{f} \tilde{\mathcal{F}}_i$ is a face
 disk of some lift of $P$.
 Since $P$ is a right-angled polyhedron, then in a neighborhood
 of $\tilde{v}$, the tesellation of $\mathbb{H}^3$
 is equivalent to the coordinate planes of $\mathbb{R}^3$.
 Since $\tilde{f} \tilde{\mathcal{F}}_i$ is a face
 disk, then in a neighborhood of $\tilde{v}$,
   $\tilde{f} \tilde{\mathcal{F}}_i$ is contained in a single
 octant.  Since $n \leq 3$ and
 $\bigcup_i \tilde{f} \tilde{\mathcal{F}}_i$ is neighborhood of
 $\tilde{f}(\tilde{v})$, then either $deg_K(v) \geq 6$ or
 $\bigcup_i \tilde{f} \tilde{\mathcal{F}}_i$
 is contained in a single octant.  However in the latter case $deg_K(v) = 3$, 
 so the face disks $\mathcal{F}_i$ are not special, for a contradiciton.
\end{proof}

\begin{lem} \label{facebound}
There is a fixed constant $c>0$ such that,
 if $f:S \rightarrow \mathbb{H}^3/\Gamma(P)$ is a cellular immersion
 with respect to a polyhedral structure $K$ on $S$, and if $f$ is
 in standard form, then $genus(S) \geq c |K^{(2)}|$.
\end{lem}

 \begin{proof}
 Since $f$ is in standard form, the 2-skeleton of $K$ has a partition
 into special face disks.  Let $K^{\prime}$ be the corresponding
 complex of $S$, and let $\mathcal{F} \in K^{\prime (2)}$.

 Let $v \in \mathcal{V}(\mathcal{F},2)$, and
 suppose $deg_{K^{\prime}}(v) \leq 3$.
 Since $deg_{\mathcal{F}}(v) = 2$, then we have
 $deg_K(v) < 2*deg_{K^{\prime}}(v) \leq 6$.
 Therefore  by Lemma \ref{deg} $deg_{K^{\prime}}(v) \geq 4$,
 for a contradiction.
 Therefore every $v \in \mathcal{V}(\mathcal{F},2)$
 has  $deg_{K^{\prime}}(v) \geq 4$.
 Since $\mathcal{F}$ is special, $|\mathcal{V}(\mathcal{F},2)| \geq 5$,
 and so

(*)\qua each $\mathcal{F}$  has at least five
 vertices $v$ with $deg_{K^{\prime}}(v) \geq 4$.

 Let $V, E$ and $F$ represent the number of vertices, edges and faces
 in the complex $K^{\prime}$.
 For $\mathcal{F} \in K^{\prime(2)}$, let $V(\mathcal{F},n)$ be the number of vertices which
 are distinct in $\mathcal{F}$ (not necessarily in $K^{\prime}$) and
 which have degree $n$ in $K^{\prime}$.
  We have:
 \begin{align*}
\chi(S) &= V - E + F\\
 &= V + (-1/2)(\sum_{v \in K^{\prime (0)}} deg(v)) + F\\
 &= -1/2(\sum_{v \in K^{\prime (0)}} deg(v)-2) +F\\
 &= -1/2 \left(
      \sum_{\mathcal{F} \in K^{\prime (2)}} 1\frac{V(\mathcal{F},3)}{3} +
      \!\!\!\!\!\sum_{\mathcal{F} \in K^{\prime (2)}} 2\frac{V(\mathcal{F},4)}{4} +
      \!\!\!\!\!\sum_{\mathcal{F} \in K^{\prime (2)}} 3\frac{V(\mathcal{F},5)}{5}
           + ...\right) + F\\
 &\leq -1/2 \left(\sum_{\mathcal{F} \in K^{\prime (2)}}
             \frac{1}{2} [V(\mathcal{F},4)+ V(\mathcal{F},5) + ...]\right)+F\\
 &\leq -1/2 \left(\frac{1}{2} 5F \right) +F, \textrm{ by (*)}\\
 &= -\frac{1}{4} F
\end{align*}

Since each face disk $\mathcal{F}$ can contain at most $|P^{(2)}|$
 faces of $P$, then $|K^{(2)}| \leq |P^{(2)}||K^{\prime (2)}|$,
 and the lemma follows.
\end{proof}

 Combining Lemma \ref{facebound} with Lemma \ref{normal}, we get the following,
 which may be of some independent interest:

\begin{pro}
 For $P$ and $M$ as in Theorem \ref{reforbs},
 there is a constant $c$ such that
 every immersed incompressible genus g surface in $M$
 can be homotoped to a union of at most $c g$ faces of $P$.
\end{pro}

We now resume the proof of the upper bound.

 By Lemma \ref{facebound},
 the total number of maximal face disks is linear in the genus,
 and so the number of edges on the boundary of these faces is linear
 in the genus also. The number of choices for the face-disks
 is then exponential in $g$,
 and the number of ways of gluing the face disks along the boundary edges
 is factorial in $g$, resulting in a factorial upper bound for
 the total number of immersions.

 We now sketch the proof of the explicit bound given in
 Theorem \ref{reforbs}.
 Let $K^{\prime}$ be a special face-disk complex for $S$.
 The above analysis gives
 $s_2(M,g) \leq e^{c_1 g}(|K^{\prime(1)}|)!$, for some constant $c_1$.
 The proof of Lemma \ref{facebound}
 shows that $\chi(S) \leq - \frac{1}{4}|K^{\prime}(2)|$,
 so $|K^{\prime (2)}| \leq 8g-8$.  Let $c(P)$ be the maximum number
 of edges in a face disk of $P$.
 Then $|K^{\prime (1)}| \leq c(P)|K^{\prime (2)}|$,
 and so $s_2(M,g) \leq e^{c_1 g}((8g-8)c(P))!$, and
 so for large $g$, $s_2(M,g) \leq e^{(8c(P)+1) g \log(g)}$.

\section{Further Questions}

Theorem \ref{main} prompts the question:

\medskip
\textbf{Question 1}\qua Does the conclusion of
 Theorem \ref{main} hold for every closed, hyperbolic 3-manifold?

\medskip
Of course, even proving that $s(g,M)$ is not constantly zero on any $M$ is
 very difficult.

 It is noteworthy that, in the examples of Theorem \ref{main},
 the factorial growth was achieved entirely by quasi-Fuchsian surfaces.
 If Thurston's virtual fibering conjecture is
 true, there should always exist in addition a certain number of geometrically
 infinite immersions.  However, the general constructions
 of $\pi_1$-injective immersions which are currently known seem to
 produce geometrically infinite examples only sporadically.

\medskip
\textbf{Question 2}\qua  Let $M$ be a closed hyperbolic 3-manifold,
 and let $i(g)$ denote the number of conjugacy classes of maximal,
 geometrically infinite surface groups of genus at most $g$
 in $\pi_1M$. Does $i(g)/s(g) \rightarrow 0$ as $g \rightarrow \infty?$
 Does $i(g)$ have sub-factorial (or even polynomial) growth?

\bibliographystyle{gtart}

\Addresses\recd
\end{document}